\documentclass[10pt,3p]{elsarticle}

\usepackage{cite}
\usepackage{amsmath,amssymb,amsfonts}
\usepackage{algorithmic}
\usepackage{graphicx}
\usepackage{textcomp}
\usepackage{xcolor}
\def\BibTeX{{\rm B\kern-.05em{\sc i\kern-.025em b}\kern-.08em
    T\kern-.1667em\lower.7ex\hbox{E}\kern-.125emX}}
\usepackage{amsmath,amssymb}
\usepackage{subfig,graphicx}
\usepackage{hyperref}
\usepackage{multirow}
\usepackage{amsthm}
\newtheorem{theorem}{Theorem}[section]
\newtheorem{lemma}[theorem]{Lemma}

\newcommand{\ff}{{\bf f}}
\newcommand{\uu}{{\bf u}}
\newcommand{\xx}{{\bf x}}
\newcommand{\mm}{{\bf m}}
\newcommand{\SO}{\text{SO}}
\newcommand{\pphi}{\boldsymbol{\phi}}

\begin{document}

\begin{frontmatter}
\title{\large\bf Closed-Form Minkowski Sums of Convex Bodies with Smooth Positively Curved Boundaries \tnoteref{t1}}

\author[1]{Sipu Ruan}
\ead{ruansp@nus.edu.sg}

\author[1,$*$]{Gregory S. Chirikjian}
\ead{mpegre@nus.edu.sg}

\address[1]{Department of Mechanical Engineering, National University of Singapore, Singapore}
\address[$*$]{Address all correspondence to this author.}
\tnotetext[t1]{This paper has been published in the journal of Computer-Aided Design. DOI: \url{https://doi.org/10.1016/j.cad.2021.103133}. The code for results in this article is available at: \url{https://github.com/ruansp/minkowski-sum-closed-form}.}
    

\begin{abstract}
This article derives closed-form parametric formulas for the Minkowski sums of convex bodies in $d$-dimensional Euclidean space with boundaries that are smooth and have all positive sectional curvatures at every point. Under these conditions, there is a unique relationship between the position of each boundary point and the surface normal. The main results are presented as two theorems. The first theorem directly parameterizes the Minkowski sums using the unit normal vector at each surface point. Although simple to express mathematically, such a parameterization is not always practical to obtain computationally. Therefore, the second theorem derives a more useful parametric closed-form expression using the gradient that is not normalized. In the special case of two ellipsoids, the proposed expressions are identical to those derived previously using geometric interpretations. In order to examine the results, numerical validations and comparisons of the Minkowski sums between two superquadric bodies are conducted. Applications to generate configuration space obstacles in motion planning problems and to improve optimization-based collision detection algorithms are introduced and demonstrated.
\end{abstract}

\begin{keyword}
Minkowski sums \sep computer-aided design \sep computational geometry 
\end{keyword}

\end{frontmatter}

\section{Introduction} \label{sec:intro}
Minkowski sums between two solid bodies have been studied for decades in the engineering and computer science literature, and have wide applications in computer-aided design and manufacturing \citep{evans1992construction}, computational geometry \citep{varadhan2004accurate}, robot motion planning \citep{lozano1990spatial}, etc. This article computes the boundary of a $d$-dimensional Minkowski sum, which is closely related to the $(d-1)$-dimensional space that describes when the two convex bodies touch each other without colliding. This space is also called the \emph{contact space} in the context of robotics \citep{behar2011fast}. As a particular example in robot motion planning, the configuration space (C-space) of a robot is constructed by computing the Minkowski sums between the robot parts at all possible orientations and the obstacles in the environment \citep{bajaj1989generation,lozano1990spatial,eckenstein2017modular}. In this case, the robot is shrunk into a point and the boundaries of obstacles are inflated, resulting in \emph{configuration-space obstacles (C-obstacles)}. Then a motion plan can be made for this point to traverse through the C-space. In addition, another class of applications includes the detection of contacts between two rigid bodies \citep{ma2018efficient}, as well as their separation distance or penetration depth \citep{lee2016Continuous}. A thorough review on this application especially for discrete bodies is provided in \citep{cox2020review}. For two convex bodies with smoothly bounded surfaces, collision detection between them can be solved by simultaneously optimizing the distance between two points, one on each body \citep{chakraborty2008proximity}. A necessary condition for this closest distance is based on the \emph{common normal} concept, which results in many efficient and elegant algorithms such as \citep{lopes2010mathematical,gonccalves2017benchmark,romer2020normal}. An equivalent result can be obtained by querying the relative position between the Minkowski sum boundary and the center point of the other body \citep{ruan2019efficient}. To achieve this result, only the point on the exact Minkowski sum boundary is optimized, as compared to simultaneously optimizing two points. It is advantageous due to the reduced numbers of variables. However, the cost of computing Minkowski sums affects the efficiency of collision detection queries. And a fast Minkowski sum computation, which is the focus of this article, can help solving collision detection problems more efficiently. The efficiency is achieved by the proposed closed-form expressions.

In general, the computations of exact Minkowski sums can be very expensive \citep{fogel2009exact}. A large number of the investigations in the literature focus on bodies that are encapsulated by polytopes, with discrete and faceted surfaces \citep{varadhan2004accurate,chirikjian2016harmonic}. They are simple to be represented and stored in a modern computer, and can characterize a large range of convex or non-convex bodies. However, this type of surface is generally not smooth and requires a large set of parameters like vertices and faces information. In this article, alternatively, we focus on a class of convex bodies whose boundary surfaces have implicit and parametric expressions. The surfaces of interest are parameterized by outward normal vectors \citep{gravesen2007surfaces,romer2020normal}, which can be either normalized or un-normalized. In other words, the \emph{inverse Gauss map} of these surfaces exists and is unique. This special property can be characterized by all sectional curvatures being positive at every point, which rules out flat or concave surface patches. This class of bodies is general and captures a wide variety of shapes with just a small number of parameters. The examples that are used to demonstrate the theories in this article include ellipsoids and superquadrics. Some basic geometric shapes can be defined easily in this class of bodies such as spheres, boxes, cylinders, etc. And with linear transformations like shearing, superquadrics can be deformed into geometries such as parallelepiped. These shapes play important roles in modeling real-world objects. They provide tight bounding volumes to capture the essential geometric properties of the object but use much fewer parameters than surface meshes. The work in this article shows that, for this wide class of convex bodies, it is possible to parameterize the boundary of their Minkowski sums in closed form. Note that the theoretical derivations throughout this article are valid for any dimension, but for real application scenarios, we only discuss the 2D and 3D cases.

Based on the two different parameterizations of the surfaces, i.e. outward normal and un-normalized gradient, two expressions for the closed-form Minkowski sums are derived, which are applicable to the bodies under general linear transformations such as rotation and shear. The novel contributions of this article are explicitly summarized as follows:
\begin{itemize}
\item Exact closed-form Minkowski sum expressions for two arbitrary convex bodies with smooth positively curved boundaries embedded in $\mathbb{R}^d$ are derived using contact-space arguments.
\item The closed-form expressions offer the flexibility of using normalized or un-normalized surface gradient parameterization.
\item The expressions under linear transformations of the two bodies are derived.
\item The applications on motion planning and collision detection problems are demonstrated.
\end{itemize}

The rest of this article is organized as follows. Section \ref{sec:literature} reviews related literature on the computational techniques for Minkowski sums. Section \ref{sec:general_mink} reviews some useful properties of surfaces that are related to the derivations in this article. Section \ref{sec:main_results} introduces the main results to compute the exact closed-form Minkowski sums between two general bodies. Section \ref{sec:demo} demonstrates the proposed expressions to the cases of ellipsoids and superquadrics. Section \ref{sec:numerial_verifications} numerically verifies the proposed method and compares the performance with the original definitions of Minkowski sums. Section \ref{sec:application} introduces two applications of the proposed method: efficient generations of configuration-space obstacles in motion planning problems; and improvements on optimization-based collision detection algorithms. Section \ref{sec:conclusion} concludes the article and points out some limitations and future work that the authors would like to further discover.

\section{Related Work} \label{sec:literature}
This section reviews related work on the computation of Minkowski sums. In general, if two bodies in $\mathbb{R}^3$ are non-convex, the complexity of computing their Minkowski sums can be as high as $O(m^3 n^3)$, where $m$ and $n$ are the features (i.e. the number of facets) of the two polytopes. And the exact complexity bounds of Minkowski sums are rigorously analyzed in \citep{fogel2009exact}. To make the calculations tractable, a large number of efficient algorithms have been proposed.

One of the most common methods for general polytopes is based on convex decomposition (either exact or approximated) since computing Minkowski sums between convex bodies can be much easier \citep{varadhan2004accurate,hachenberger2009exact}. For example, a direct usage of the definition of Minkowski sums and convex hull algorithms provides a straight-forward way to compute Minkowski sums between two convex discrete bodies in any dimension \citep{lozano1990spatial}. The idea is to add all pairs of points on the surface of each body, resulting in a superset of the actual Minkowski sums boundary. Then, the boundary can be computed by taking the convex hull of this superset. Section \ref{sec:general_mink:def} provides a more detailed mathematical definition of this method and Fig. \ref{fig:mink:demo_mink_sum:def} illustrates the idea. This method is chosen as a baseline when conducting benchmark simulations in this article. Its performance depends highly on the convex hull algorithms, such as Graham's \citep{graham1983finding}, QuickHull \citep{barber1996quickhull}, etc. But comparisons among different implementations of convex hulls are beyond the scope of this article. The core idea behind decomposition-based methods is that the union of Minkowski sums is the Minkowski sums of the union of bodies. Therefore, for this type of methods, the efficiency of convex decomposition and union operation affects the overall performance of the algorithms. Although the focus of this article is convex bodies, it is still possible to extend the capability of our work to deal with non-convex bodies using convex decomposition methods.

A unified framework based on \emph{slope diagram} has been presented for both convex and non-convex bodies in both 2D and 3D spaces \citep{ghosh1993unified}. The idea is closely related to the \emph{Gauss map}, which maps a surface point into its normal vector. This framework is able to deal with polygons/polyhedra or bodies with smooth bounding curves. The slope diagrams of both bodies are computed and merged together via stereographic projection technique. Then the Minkowski sums can be computed by extracting the boundary from the merged diagram. Concretely for 2D convex polygons, an extreme point on the Minkowski sum boundary at one direction is the sum of extreme points on the two bodies along that direction. Therefore, by adding vertices of each body at sorted directions according to angles between edges, this algorithm achieves linear time complexity. This algorithm is involved in the benchmark studies of this article for 2D cases, which is denoted as \emph{edge sorting} \citep{mark2008computational}. Works on improving this framework have also been done. For example, the expensive stereographic projection step for merging the slope diagram can be more efficient by using vector operations in the case of convex polyhedra \citep{wu2003improvements}. More improvements are made using the \emph{contributing vertices} concept \citep{barki2009contributing}. It directly finds the extreme vertex on the other body along a specific direction, therefore avoids explicitly computing slope diagram and stereographic projection. Other similar work for convex polyhedra use variants of Gauss map such as \citep{fogel2007exact}. It proposes the idea of \emph{cubical Gaussian map}, which is a projection of Gauss map onto the smallest axis-aligned unit bounding cube that contains it. Like the original Gauss map, a point on the unit cube is related to the unit normal vector of the body facet. This article applies the Gauss map concept extensively. One of the sufficient conditions for the existence of the proposed closed-form expression is the uniqueness of the inverse Gauss map. Different from searching the corresponding vertices and facets, by using the inverse Gauss map, the extreme vertex of one body can be explicitly and directly expressed by the parameter of the other body in closed form.

Minkowski sums of two bodies can also be formulated using convolution. The Minkowski sums of two solid bodies is the support of the convolution of their indicator functions \citep{chirikjian2016harmonic}. The core idea of this type of methods is to match the normal vectors of the two body surfaces. And if two bodies are convex, the convolution is exactly the same as the Minkowski sum boundary. In \citep{kavraki1995computation}, the convolution of two bodies is computed using fast Fourier transforms. In \citep{lien2009simple}, all the possible pairs of valid features are computed in a brute-force way, then those that are not on the boundary surface are trimmed out using collision detection methods. Recent work of \citep{baram2018exact} further extends the procedure using reduced convolution via some filters, and applies to bodies with holes inside. Convolution-based methods are also widely applied between surfaces or curves with algebraic expressions. Earlier work such as \citep{lee1998polynomial} propose several algorithms to approximate the exact convolution between two planar algebraic curves using polynomial/rational curves. In \citep{muhlthaler2003computing}, Minkowski sums are further studied for bodies bounded by surface patches with linear normal vector field. In \citep{peternell2007minkowski}, the smooth or piecewise smooth surfaces, represented by dense point clouds and triangulated facets, are studied. Recent work of \citep{mizrahi2017minkowski} computes Minkowski sums for parametric B-spline surfaces. This article uses the similar spirit of matching normal vectors at the the two surface points. Compared to these work in literature, this article studies the bodies that can be explicitly and uniquely parameterized by surface normal in closed form, such as ellipsoids and superquadrics. We further relax the normalization by matching the surface gradient of the two bodies at the touching point in closed form. Also, this work shows the closed-form expression under linear transformations such as rotation and shearing deformation.

The most relevant work with this article is \citep{yan2015closed}, where closed-form expression of Minkowski sums between two ellipsoids has been proposed. The method used there was based on geometric interpretations using affine transformation of ellipsoids. By transforming the ellipsoid into a sphere, the Minkowski sums in the affine space can be computed as an offset surface in closed form. The exact Minkowski sums boundary can then be obtained by transforming the whole space back. The expressions are exact, in closed-form and valid in any dimension. Extension to this work has also been proposed when one of the bodies is a general convex smooth one \citep{ruan2018path}. However, both these work require at least one body to be an ellipsoid, which limits the ability to extend into more complex shapes. This article also develops exact and closed-form Minkowski sums, but allows both two bodies to be more general. The resulting expressions in the cases that involve ellipsoid(s) are identical to \citep{yan2015closed,ruan2018path}, which are just special cases of this work. Moreover, this work derives these expressions in totally different ways.

\section{Minkowski Sums and Geometric Properties of Surfaces} \label{sec:general_mink}

This section reviews definitions and properties of Minkowski sums and surfaces that are necessary for the derivations throughout this article. 

\subsection{Definition and computation of Minkowski sums} \label{sec:general_mink:def}
Given two solid bodies $B_1$ and $B_2$ in $d$-dimensional Euclidean space, their Minkowski sum is defined as
\begin{equation}
B_1 \oplus B_2 \,\doteq\, \{\xx + {\bf y} \,|\, \xx \in B_1\,,\, {\bf y} \in B_2\}\,.
\label{eq:mink_sum_def}
\end{equation}
In this article, we specifically compute the boundary of Minkowski sums between $B_1$ and $-B_2$, i.e. $\partial [B_1 \oplus (-B_2)]$, where
\begin{equation}
B_1 \oplus (-B_2) \,\doteq\, \{\xx - {\bf y} \,|\, \xx \in B_1\,,\, {\bf y} \in B_2\}\,.
\label{eq:our_result_def}
\end{equation}
Geometrically, when placing $B_2$ at this boundary surface, $B_1$ and $B_2$ touch each other (without penetrating) outside. In the context of robotics, this boundary is also called the \emph{contact space}, which divide the space of collision and separation. Note that $-B_2$ is the reflection of the original body $B_2$ with respect to its center. And if $B_2$ is centrally symmetric to itself (like an ellipsoid), Eqs. \eqref{eq:mink_sum_def} and \eqref{eq:our_result_def} are equivalent, i.e. $B_1 \oplus B_2 = B_1 \oplus (-B_2)$. Once this boundary is known, $B_1 \oplus (-B_2)$ can also be parameterized trivially.

If a set of points is sampled on the boundary of each surface, the Minkowski sums between these two bodies can be approximated as a discrete surface. Since both boundary surfaces are assumed to be convex, a direct way to compute the approximated Minkowski sums is to add all the points on different surfaces and take the \emph{convex hull} of the result \citep{fogel2007exact}. This refers to a direct method using the definition, i.e. Eq. \eqref{eq:our_result_def}, of the Minkowski sums, whose performance is compared with our closed-form method in Sec. \ref{sec:numerial_verifications}.

\begin{figure*}
\centering
\subfloat[Computational process via the original Minkowski sums definition.]{\includegraphics[scale=0.3, trim=45 150 60 100, clip]{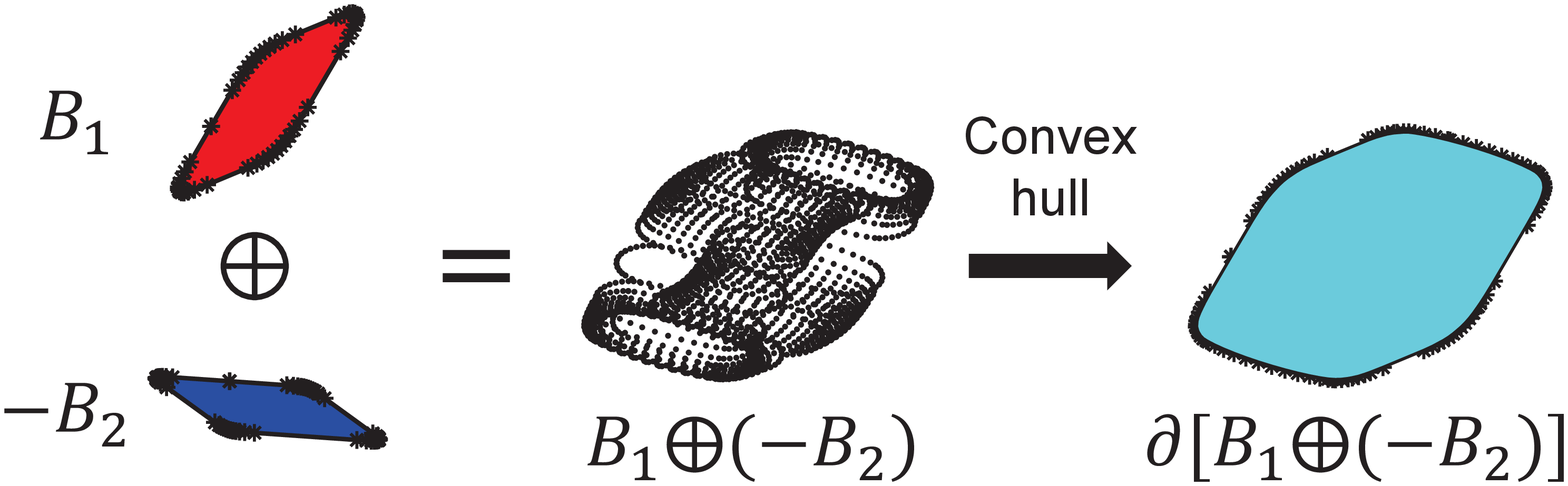}
\label{fig:mink:demo_mink_sum:def}}
~
\subfloat[Computational process of the closed-form solution.]{\includegraphics[scale=0.3, trim=150 120 150 100, clip]{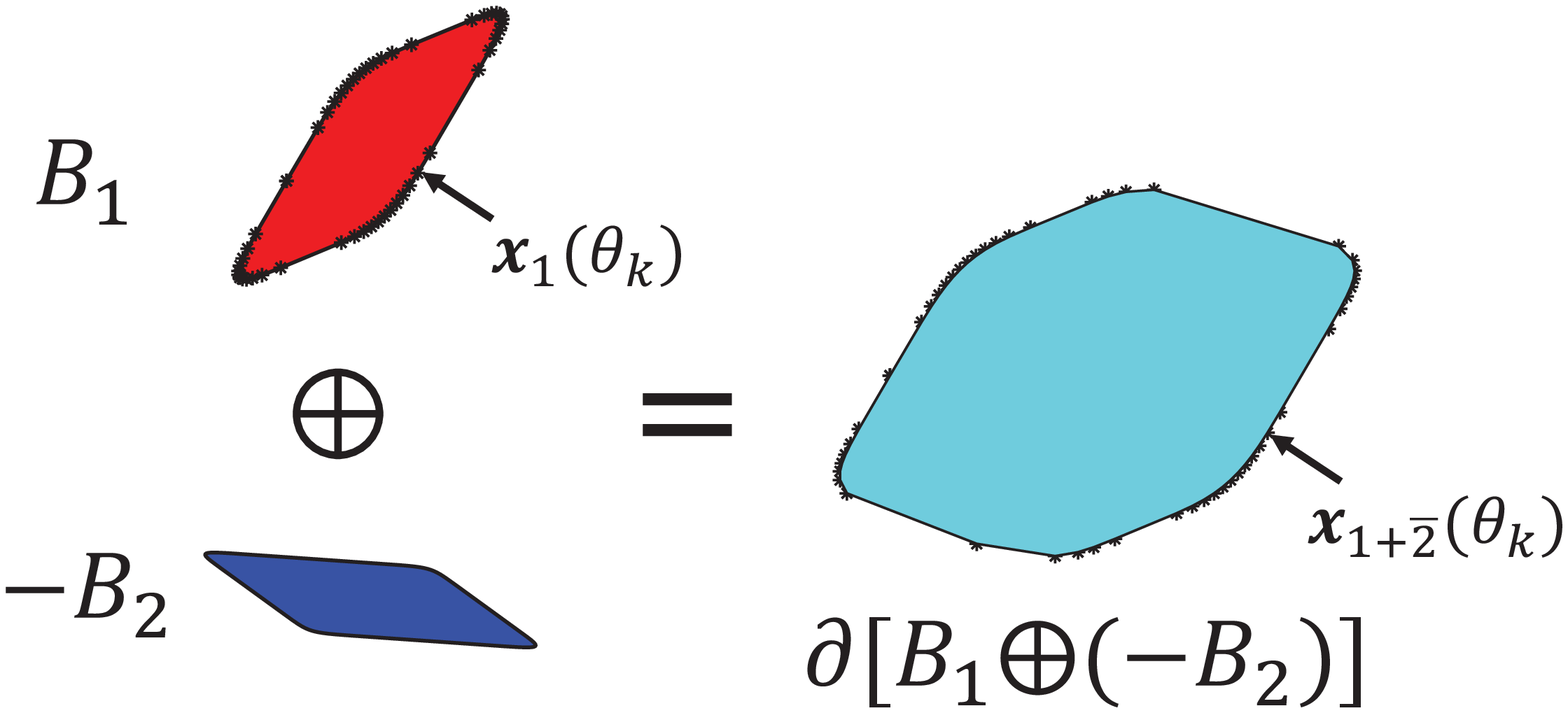}}
\caption{Demonstration of the process to compute Minkowski sums using the original definition and our closed-form solutions. In (a), the surfaces enclosing both bodies are sampled into point sets; then all the pairs of points in different surfaces are added to obtain the Minkowski sums, i.e. $B_1 \oplus (-B_2)$; the boundary of Minkowski sums, i.e. $\partial [B_1 \oplus (-B_2)]$, is then computed via convex hull operation. In (b), only one of the body (i.e. $B_1$) surfaces is sampled, and the parameters (i.e. $\theta_k$) of the points on Minkowski sums boundary is the same with those on the sampled body; in other words, $\xx_{1+\overline{2}}(\theta_k)$ is computed based on $\xx_1(\theta_k)$ in closed-form. Here, $\xx_{1+\overline{2}}$ denotes the boundary point on $\partial [B_1 \oplus (-B_2)]$, with $\overline{2}$ representing the reflected body, i.e. $-B_2$. And $\theta_k \in [-\pi, \pi]$ is the sampled angle that parameterizes the boundary of a unit circle, i.e. $[\cos \theta_k, \sin \theta_k]^T$.}
\label{fig:mink:demo_mink_sum}
\end{figure*}

Figure \ref{fig:mink:demo_mink_sum} shows the process of computing Minkowski sums using the original definition as in Eq. \eqref{eq:our_result_def} as well as our closed-form solution. The demonstrations in the 2D case are shown in the figure, and the boundaries of both bodies as well as their Minkowski sums are 1D curves. For both methods, the points on the bounding surfaces are sampled based on a set of angles $\theta_k \in [-\pi,\pi]$ that parameterize an 1D unit circle.

\subsection{Some geometric properties of surfaces}
Suppose that the boundary of convex closed body $B_i$ is defined by the parametric equation
\begin{equation}
\xx = {\bf f}({\bf u}) \,\in\, \mathbb{R}^d
\label{eq:param_equation}
\end{equation}
where the unit vector ${\bf u} \in \mathbb{S}^{d-1}$ is in turn parameterized by $d-1$ angles. A coordinate system can always be chosen such that the body $B_i$ itself can be centered and parameterized as $\xx = r \, {\bf f}({\bf u})$ where $r \in [0,1]$. Then the origin ${\bf 0} \in B_i$. Here ${\bf f}:\mathbb{S}^{d-1} \,\rightarrow\, \mathbb{R}^d$ can be thought of as an embedding of the sphere which is deformed while ensuring the convexity of $B_i$. Consequently, the line segment connecting any two points in $\partial B_i$ is fully contained in $B_i$.

Assume that the corresponding implicit equation of $B_i$ exists and can be written as
\begin{equation}
\Psi(\xx) \,=\, \Psi({\bf f}({\bf u})) \,=\, 1 \,.
\label{eq:implicit_equation}
\end{equation}

Thinking of ${\bf u} = {\bf u}(\phi_1,...,\phi_{d-1})$, the tangent vectors can be computed from Eq. \eqref{eq:param_equation} as ${\bf t_i} = \frac{\partial \xx}{\partial \phi_i}$, which, in general, is not orthonormal. Moreover, the outward pointing unit normal to the surface at the same point $\xx$ can be calculated as
\begin{equation}
{\bf n}({\bf u}) \,\doteq\, \left. \frac{(\nabla_{x} \Psi)(\xx)}{\|(\nabla_{x} \Psi)(\xx)\|} \right|_{\xx = {\bf f}({\bf u})}  \,.
\label{eq:outward_normal}
\end{equation}
The \emph{Gauss map} assigns to each point ${\bf f}({\bf u}) \in \partial B$ its normal ${\bf n}({\bf u}) \in \mathbb{S}^{d-1}$. And the function ${\bf n}: \mathbb{S}^{d-1} \,\rightarrow\, \mathbb{S}^{d-1}$ is closely related. In this article we will make extensive use of Eq. \eqref{eq:outward_normal}, as well as the un-normalized gradient
\begin{equation}
\mm({\bf u}) \,\doteq\, \left. (\nabla_{x} \Psi)(\xx) \right|_{\xx = {\bf f}({\bf u})}  \,.
\label{eq:unnormalized_normal}
\end{equation}

The Gaussian curvature is the Jacobian determinant of the Gauss map. Consequently,
if $\partial B$ has all positive sectional curvature everywhere, then the relationship
${\bf n} = {\bf n}({\bf u})$ will be invertible as ${\bf u} = {\bf n}^{-1}({\bf u}) \doteq {\bf u}({\bf n})$. Then it becomes possible, at least in principle, to re-parameterize positions on the surface using the outward normal as
\begin{equation}
\xx = \tilde{{\bf f}}({\bf n}) \,,
\label{eq:param_equation_normal}
\end{equation}
where $\tilde{{\bf f}}({\bf n}) = {\bf f}({\bf u}({\bf n}))$. As with ${\bf u}$, it is possible to use angles to parameterize ${\bf n}$. When such a parameterization is used,
\begin{equation}
\left.\frac{(\nabla_{x} \Psi)(\xx)}{\|(\nabla_{x} \Psi)(\xx)\|} \right|_{\xx = \tilde{{\bf f}}({\bf n})} \,=\, {\bf n} \,.
\label{eq:normal_reparam}
\end{equation}

\section{Main Results} \label{sec:main_results}

\begin{figure*}
\centering
\subfloat[Visualization of Theorem \ref{thm:normal} \label{fig:demo_thm_1}]{\includegraphics[scale=0.6, trim=60 20 60 20, clip]{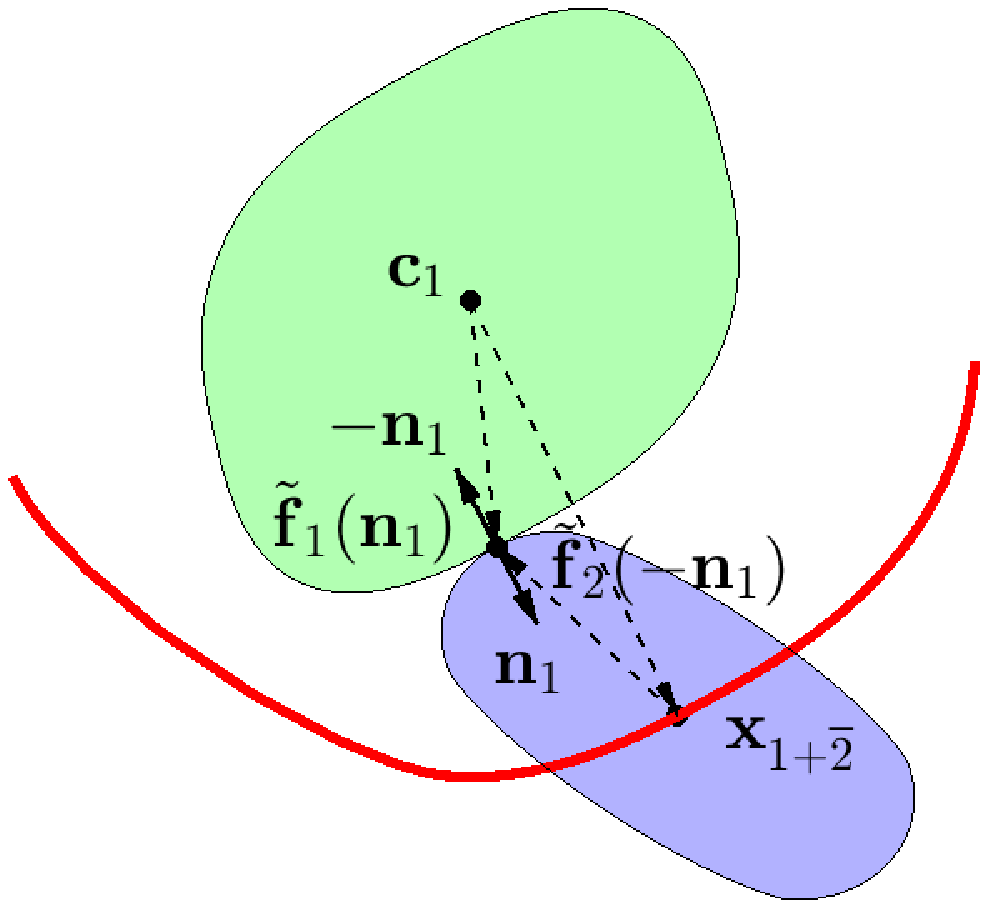}}
~
\subfloat[Visualization of Theorem \ref{thm:gradient} \label{fig:demo_thm_2}]{\includegraphics[scale=0.6, trim=60 20 60 20, clip]{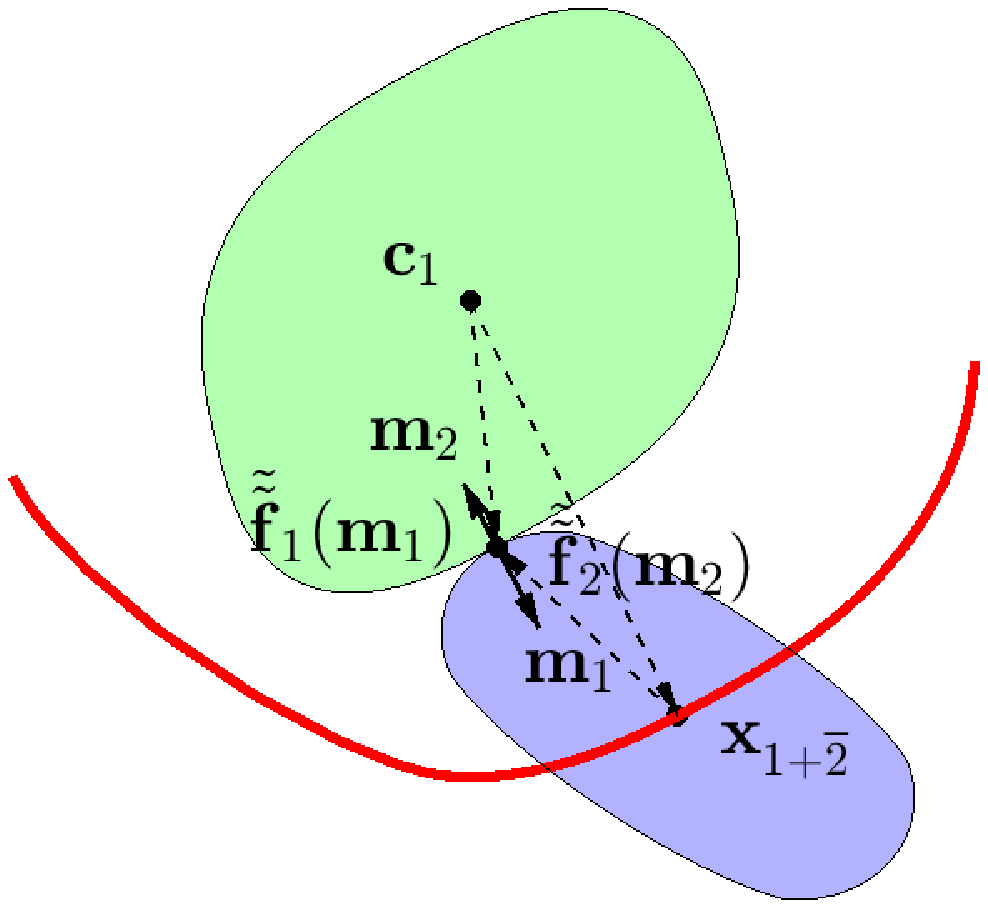}}
\caption{Demonstration of the main results. The two bodies touch each other at only one point, which has different parameterizations. (a) Theorem \ref{thm:normal}: Two bodies are parameterized by the unit outward normal vectors, i.e. $\tilde{{\bf f}}_1({\bf n}_1)$ and $\tilde{{\bf f}}_2({-\bf n}_1)$. The normal vectors at the contact point are anti-parallel to each other. (b) Theorem \ref{thm:gradient}: Two bodies are parameterized by the un-normalized gradients ($\mm_1$ and $\mm_2$), where $\mm_2 = -\frac{\Phi(\mm_1)}{\| \mm_1 \|} \mm_1$. They are anti-parallel but have different magnitudes. The green and blue patched bodies are $B_1$ and $B_2$ respectively. ${\bf c}_1$ is the center of $B_1$, and $B_2$ is placed at one of the Minkowski sums boundary point $\xx_{1+\overline{2}}$. The red curve is a section of the computed closed-form Minkowski sum boundary. }
\label{fig:demo_thms}
\end{figure*}

In this section, the main results of the article are presented. We first propose two possible expressions for closed-form Minkowski sums, along with the proofs. Both of these theorems assume two bodies are in their canonical forms. Then, we show the expression when linear transformations are applied to the two bodies. The expressions are demonstrated in the 2D case in Fig. \ref{fig:demo_thms}.

\subsection{Closed-form Minkowski sums parameterized by outward unit normal line space}

\begin{theorem}
If the boundaries $\partial B_1$ and  $\partial B_2$ of convex bodies $B_1$ and $B_2$ each have closed form parametric expressions of the form in Eq. \eqref{eq:param_equation} and their Gauss maps can be obtained in closed form resulting in parametric equations $\xx_i = \tilde{{\bf f}}_i({\bf n}_i)$, then their contact space $\partial[B_1 \oplus (-B_2)]$ can be parameterized as
\begin{equation}
\boxed{
\xx_{1+\overline{2}}({\bf n}_1) = \tilde{{\bf f}}_1({\bf n}_1) - \tilde{{\bf f}}_2(-{\bf n}_1) } \,,
\label{eq:mink_sum_norm}
\end{equation}
where ${\bf n}_1 \in \mathbb{S}^{d-1}$ can be parameterized using spherical coordinates which in turn parameterize $\xx_{1+\overline{2}}$ \footnote{Throughout this article, the boundary point of $B_1 \oplus B_2$ is denoted as $\xx_{1+2}$, while the boundary point of $B_1 \oplus (-B_2)$ is denoted as $\xx_{1+\overline{2}}$. If a body is centrally symmetric, like an ellipsoid, $\xx_{1+2} = \xx_{1+\overline{2}}$.}.
\label{thm:normal}
\end{theorem}

\begin{proof}
Hold $B_1$ fixed, and allow $B_2$ to move such that it kisses $B_1$ at any point on its surface. Record the surface that the origin of the reference frame attached to $B_2$ traces out as $B_2$ undergoes all translational motions for which $\partial B_1 \cap \partial B_2$ consists of a single point. When this occurs, the normals of both surfaces must pass through the point of contact, lie along a common line, and have opposite sense. The common point of contact between the two bodies will appear as $\tilde{{\bf f}}_1({\bf n}_1)$ in the coordinate system of $B_1$ (i.e., the world frame), and as $\tilde{{\bf f}}_2(-{\bf n}_2)$ is the moving coordinate system attached to $B_2$. The latter is described in the coordinate system in the world frame as the relative position vector $-\tilde{{\bf f}}_2(-{\bf n}_2)$. Adding both gives the result.
\end{proof}

Obviously, given this contact space, replacing $B_2$ with $-B_2$ then gives the Minkowski sum boundary $\partial (B_1 \oplus B_2)$ as $\xx_{1+2}({\bf n}_1) = \tilde{{\bf f}}_1({\bf n}_1) + \tilde{{\bf f}}_2({\bf n}_1)$. Figure \ref{fig:demo_thm_1} demonstrates the relationships of the outward normal vectors at each body boundary (i.e. ${\bf n}_1$) as well as different terms in parameterizing $\partial[B_1 \oplus (-B_2)]$ as in Eq. \eqref{eq:mink_sum_norm}. Though this can be done in principle, in practice it is often difficult to obtain in closed form expression for $\tilde{{\bf f}}_i({\bf n}_i)$.

In such cases, it can be easier to parameterize $\xx$ in terms of its un-normalized gradient as
\begin{equation}
\xx = \tilde{\tilde{{\bf f}}}(\mm) \,,
\label{eq:param_equation_gradient}
\end{equation}
where $\mm$ is obtained from Eq. \eqref{eq:unnormalized_normal}. Then it is possible to equate
$ {\bf f}({\bf u}) \,=\, \tilde{{\bf f}}({\bf n}) \,=\,  \tilde{\tilde{{\bf f}}}(\mm) $,
because Eqs. \eqref{eq:param_equation}, \eqref{eq:param_equation_normal} and \eqref{eq:param_equation_gradient} are different parameterizations of the same boundary, i.e. $\partial B$. Once $\tilde{\tilde{{\bf f}}}(\mm)$ is obtained analytically, an explicit surface parameterization can be constructed by setting
\begin{equation}
\mm(\phi) = \tilde{\tilde{{\bf f}}}^{-1}({\bf f}({\bf u}(\phi))) \,.
\label{eq:m_inv}
\end{equation}
In the case of Eq. \eqref{eq:param_equation_normal}, there is greater choice: Either the analogous computation can be done as
\begin{equation}
{\bf n}(\phi) = \tilde{{\bf f}}^{-1}({\bf f}({\bf u}(\phi))) \,,
\label{eq:n_inv}
\end{equation}
or we can simply choose to parameterize ${\bf n}$ in the same way as ${\bf u}(\phi)$. That is, even though ${\bf n}$ is not ${\bf u}$, in the final step ${\bf u}$ is no longer part of the description and we can choose to parameterize the surface by letting ${\bf n} = {\bf u}(\phi)$. Hence these are two variants on this parameterization. The same choice does not exist for parameterizing $\mm(\phi)$ because it is not a unit vector, and Eq. \eqref{eq:m_inv} is used. Note that by definition 
\begin{equation}
 \left. (\nabla_{x} \Psi)(\xx) \right|_{\xx = \tilde{\tilde{{\bf f}}}(\mm)} \,=\, \mm \,.
\label{eq:gradient_reparam}
\end{equation}
And this leads to our second main result.

\subsection{Closed-form Minkowski sums parameterized by un-normalized gradient}
Before presenting the result, the following lemma provides a sufficient condition for the existence of the closed-form expression of $\mm_2$ with respect to $\mm_1$.

\begin{lemma}
Suppose that $B_1$ and $B_2$ touch each other at a point, and $\mm_1, \mm_2$ are the gradient vectors at the contact point of the two bodies respectively. The closed-form expression of $\mm_2$ with respect to $\mm_1$ exists if the following two conditions are both satisfied: 
\\
(1) $\mm_2$ can be expressed in closed-form as a function of the unit hypersphere $\uu_2$, i.e. $\mm_2 = {\bf g}_2(\uu_2)$; \\
(2) ${\bf g}_2(k \uu_2) = \iota(k) \, {\bf g}_2(\uu_2)$, where $k$ is a constant scalar and $\iota(k)$ is a constant function of $k$.
\label{lemma:gradient:suff_cond}
\end{lemma}

\begin{proof}
At the contact point, the outward normal vectors of the two surfaces are anti-parallel, with the relationship being
\begin{equation}
\frac{\mm_2}{\|\mm_2\|} = - \frac{\mm_1}{\|\mm_1\|}.
\label{eq:normal_equation}
\end{equation}

If condition (1) is satisfied, i.e. $\mm_2 = {\bf g}_2(\uu_2)$, then Eq. \eqref{eq:normal_equation} becomes
\begin{equation}
\frac{ {\bf g}_2(\uu_2) }{\|\mm_2\|} = - \frac{\mm_1}{\|\mm_1\|}.
\end{equation}
This means that ${\bf g}_2(\uu_2)$ is proportional to $-\mm_1$. If condition (2) is satisfied, i.e. ${\bf g}_2(k \uu_2) = \iota(k) \, {\bf g}_2(\uu_2)$, then $\uu_2$ is proportion to ${\bf g}_2^{-1}( -\mm_1 )$. Therefore, using the fact that $\uu_2$ is a unit vector, we get
\begin{equation}
\uu_2 = \frac{ {\bf g}_2^{-1}(\mm_1) }{ \| {\bf g}_2^{-1}(\mm_1) \| }.
\end{equation}
Then,
\begin{equation}
\| \mm_2 \| = \| {\bf g}_2(\uu_2) \| = \left\| {\bf g}_2 \left( \frac{ {\bf g}_2^{-1}(\mm_1) }{ \| {\bf g}_2^{-1}(\mm_1) \| } \right) \right\| \doteq \Phi(\mm_1) \,,
\label{eq:Phi_m1_general}
\end{equation}
which concludes the proof.
\end{proof}

\begin{theorem}
Suppose the boundaries $\partial B_1$ and  $\partial B_2$ of convex bodies $B_1$ and $B_2$ each have closed form parametric expressions of the form in Eq. \eqref{eq:param_equation_gradient} and if Lemma \ref{lemma:gradient:suff_cond} is satisfied, i.e. $\|\mm_2\|=\Phi(\mm_1)$ can be written in closed form, then the contact space $\partial[B_1 \oplus (-B_2)]$ can be parameterized in closed form as
\begin{equation}
\boxed{
\xx_{1+\overline{2}}(\mm_1) = \tilde{\tilde{{\bf f}}}_1(\mm_1) - \tilde{\tilde{{\bf f}}}_2\left(-\frac{\Phi(\mm_1)}{\|\mm_1\|} \mm_1\right)
} \,,
\label{eq:mink_sum_non_norm}
\end{equation}
where $\mm_1 \in \mathbb{R}^{d}$ can be parameterized using the spherical coordinates of the original parametric expression for $B_1$.
\label{thm:gradient}
\end{theorem}

\begin{proof}
Using the same reasoning as in Theorem \ref{thm:normal}, for the bodies to kiss at a point their gradients must lie along the common normal line passing through the kissing point, and have opposite sense. Since $\mm_i = \|\mm_i\| {\bf n}_i$ ($i = 1,2$), the appropriate match of gradients that will cause a point $\xx_2(\mm_2)$ in the frame attached to $B_2$ to kiss $B_1$ at the point $\xx_1(\mm_1)$ in the frame attached to $B_1$, can be screened for by computing
$$ \int_{\mathbb{S}^{d-1}} \delta({\bf n}_1,{\bf n}_2) \tilde{\tilde{{\bf f}}}_2(-\|\mm_2\| {\bf n}_2) \, d{\bf n}_2
\,=\, \tilde{\tilde{{\bf f}}}_2(-\|\mm_2\| {\bf n}_1) \,. $$
This can then be written as
$$  \tilde{\tilde{{\bf f}}}_2\left(-\frac{\|\mm_2\|}{\|\mm_1\|} \mm_1\right) \,=\,
\tilde{\tilde{{\bf f}}}_2\left(-\frac{\Phi(\mm_1)}{\|\mm_1\|} \mm_1\right) \,. $$
The existence of this result can be achieved sufficiently using Lemma \ref{lemma:gradient:suff_cond}. Then, subtracting gives Eq. \eqref{eq:mink_sum_non_norm}.
\end{proof}

Consequently, $\xx_{1+2} = \tilde{\tilde{{\bf f}}}_1(\mm_1) + \tilde{\tilde{{\bf f}}}_2\left(\frac{\Phi(\mm_1)}{\|\mm_1\|} \mm_1\right)$. From these theorems, it is easy to observe that when $B_2 = c \, B_1$, $\xx_{1+2} = (1+c) \xx_1$. Figure \ref{fig:demo_thm_2} shows the relationships of each term in the expression derived in Eq. \eqref{eq:mink_sum_non_norm}.

\subsection{General closed-form Minkowski sum expressions with linear transformations}

Under invertible linear transformations from $\mathbb{R}^d$  defined by
$\xx' = M \xx$, tangents transform linearly as ${\bf t}_i \,\rightarrow\, M {\bf t}_i$. Alternatively, gradients transform as
\begin{equation}
\left.(\nabla_{x'} \Psi)(M^{-1} \xx')\right|_{\xx' = M \xx} \,=\, M^{-T} (\nabla_{x} \Psi)(\xx) \,.
\label{eq:gradient_trans}
\end{equation}
This can be written as
\begin{equation}
\mm' \,=\, M^{-T} \mm \,,
\label{eq:m_prime}
\end{equation}
where $M^{-T} \doteq (M^T)^{-1} = (M^{-1})^T$. This can be used together with $\tilde{\tilde{{\bf f}}}_i(\mm_i)$ to get the Minkowski sum of linearly transformed parametric surfaces, by observing that the result of transforming by $M$ is
$$ 
\tilde{\tilde{{\bf f}}}'_i(\mm'_i) \,=\, M_i \tilde{\tilde{{\bf f}}}_i(\mm_i)  \,=\, M_i \tilde{\tilde{{\bf f}}}_i(M_{i}^{T}\mm'_i) \,. 
$$

Then substituting either of these into Eq. \eqref{eq:mink_sum_non_norm} and expanding gives 
\begin{equation}
\boxed{\,
\begin{aligned}
\xx'_{1+\overline{2}}(\mm'_1) = M_1 \tilde{\tilde{{\bf f}}}_1(M_{1}^{T} \mm'_1) - M_2 \tilde{\tilde{{\bf f}}}_2 \left( -\frac{\Phi( M_{2}^{T} \mm'_1)}{\|M_{2}^{T} \mm'_1\|} M_{2}^{T} \mm'_1 \right)
\end{aligned}
 \,} \,,
\label{eq:mink_sum_transform_1}
\end{equation}
and
\begin{equation}
\boxed{\,
\begin{aligned}
\xx'_{1+\overline{2}}(\mm_1) = M_1 \tilde{\tilde{{\bf f}}}_1(\mm_1) - M_2 \tilde{\tilde{{\bf f}}}_2 \left( -\frac{\Phi( M_{2}^{T} M_{1}^{-T} \mm_1)}{\|M_{2}^{T} M_{1}^{-T} \mm_1\|} M_{2}^{T} M_{1}^{-T} \mm_1 \right)
\end{aligned}
\,} \,.
\label{eq:mink_sum_transform_2}
\end{equation}

\section{Demonstrations with Ellipsoids and Superquadrics} \label{sec:demo}

This section demonstrates the proposed closed-form Minkowski sums between two ellipsoids and two superquadrics. Only some essential expressions are derived and shown here. For more detailed derivations of intermediate steps, please refer to the supplementary materials.

\subsection{Demonstrations with ellipsoids} \label{sec:demo:ellipsoid}

Two theorems are demonstrated with ellipsoids. The results are compared with the ones coming from the previous geometric interpretations, which are the same. But using the derivations from this article, the expressions are symmetrical and reflects the commutative property of Minkowski sums.

The surface of an ellipsoid can be parameterized as
\begin{equation}
{\bf f}({\bf u}) = A {\bf u} \,\in\, \mathbb{R}^d
\label{eq:ellipsoid_param}
\end{equation}
where $A= A^T \in \mathbb{R}^{d \times d}$ is a positive definite matrix, and ${\bf u} \in \mathbb{S}^{d-1}$, which can be parameterized by spherical angles $\phi = (\phi_1, ..., \phi_{d-1})$. The corresponding implicit equation is
\begin{equation}
\Psi(\xx) \,=\, 1 \,\,\, {\rm where} \,\,\,
\Psi(\xx) \,\doteq\, \xx^T A^{-2} \xx \,.
\label{eq:ellipsoid_implicit}
\end{equation}

\subsubsection{Demonstration of Theorem \ref{thm:normal}}

By writing ${\bf u}$ as a function of ${\bf n}$, we have
$$ \tilde{{\bf f}}({\bf n})  \,=\, A \frac{A {\bf n}}{\|A {\bf n}\|} \,. $$
Then, from Theorem \ref{thm:normal},

\begin{equation}
\xx_{1+\overline{2}}({\bf n}_1) \,=\,  \frac{A_1^2 {\bf n}_1}{\|A_1 {\bf n}_1\|} \,+\, \frac{A_2^2 {\bf n}_1}{\|A_2 {\bf n}_1\|}  \,.
\label{eq:ellipsoid_mink_sum_thm_1}
\end{equation}
This is a variant of the formular given in \citep{yan2015closed} and the same as the formula derived in \citep{chirikjian2021applications} \footnote{Since ellipsoid is centrally symmetric, the expression derived here is equivalent with the one in \citep{chirikjian2021applications}, i.e. $\xx_{1+\overline{2}} = \xx_{1+2}$.}. The unit vector ${\bf n}_1 \in \mathbb{S}^{d-1}$ can be parameterized with spherical angles $\phi = (\phi_1, ..., \phi_{d-1})$ just like ${\bf u}(\phi)$, though it should not be confused with it.

\subsubsection{Demonstration of Theorem \ref{thm:gradient}}

An alternative would be to parameterize ellipsoids according to Theorem \ref{thm:gradient}. For two ellipsoids, i.e. $i = 1,2$, the surface function with respect to the gradient can be computed as
\begin{equation}
\tilde{\tilde{{\bf f}}}_i(\mm_i) = \frac{1}{2} A_{i}^{2} \mm_i \,.
\label{eq:ellipsoid_param_m}
\end{equation}

$B_1$ and $B_2$ will touch at a surface point when 
$$  \frac{\mm_1}{\|\mm_1\|} \,=\, -\frac{\mm_2}{\|\mm_2\|} \,. $$
Substituting the gradient expression $ \mm_i \,=\, \tilde{\tilde{{\bf f}}}_{i}^{-1}({\bf f}({\bf u}_i)) = 2 A_{i}^{-1} {\bf u}_i $
into the touching condition gives
$$ \frac{A_{1}^{-1} {\bf u}_1}{\|A_{1}^{-1} {\bf u}_1\|} \,=\, - \frac{A_{2}^{-1} {\bf u}_2}{\|A_{2}^{-1} {\bf u}_2\|} \,. $$
This equation says that ${\bf u}_2$ must be proportional to $-A_2 A_{1}^{-1} {\bf u}_1$, and the constraint that ${\bf u}_2$ is a unit vector then specifies that
$$ {\bf u}_2 \,=\, -\frac{A_2 A_{1}^{-1} {\bf u}_1}{\|A_2 A_{1}^{-1} {\bf u}_1\|}\,. $$
Then, the expression of $\Phi(\mm_1)$ can be obtained as
$$ \Phi(\mm_1) \,\doteq\, \|\mm_2\| \,=\, \frac{2 \|\mm_1\|}{\|A_2 \mm_1 \|} \,. $$
Therefore, Theorem \ref{thm:gradient} states that 
\begin{equation}
\xx_{1+\overline{2}}(\mm) \,=\, \frac{1}{2} A_1^2 \mm_1 + A_2^2  \frac{\mm_1}{\|A_2 \mm_1 \|} \,.
\label{eq:ellipsoid_mink_sum_closed_m}
\end{equation}

Substituting $\mm_1 = 2 A_{1}^{-1} {\bf u}_1$ with ${\bf u}_1$ parameterized as ${\bf u}_1(\phi)$ gives 
\begin{equation}
\xx_{1+\overline{2}}(\phi) \,=\, A_1 {\bf u}(\phi) +  \frac{A_2^2 A_{1}^{-1} {\bf u}(\phi)}{\|A_2 A_{1}^{-1} {\bf u}(\phi) \|} \,.
\label{eq:ellipsoid_mink_sum_phi}
\end{equation}
This is the same as the closed-form formula derived in a completely different way in \citep{yan2015closed}.

\subsubsection{Demonstration of linear transformations}
Ellipsoids can be viewed as linearly transformed versions of spheres. For a sphere, ${\bf f}({\bf u}) = {\bf u}$
and since ${\bf n} = \mm = {\bf u}$, we have
$\tilde{{\bf f}}({\bf n}) = {\bf n}$ and $\tilde{\tilde{{\bf f}}}(\mm) \,=\, \mm$. Viewing the Minkowski sum of two ellipsoids as that of two transformed spheres, where $M_i = A_{i}$, then $\mm'_i = A_{i}^{-1} {\bf u}$ and $\mm_i = {\bf u}$ and Eqs. \eqref{eq:mink_sum_transform_1} and \eqref{eq:mink_sum_transform_2} both reduce to (\ref{eq:ellipsoid_mink_sum_closed_m}).

\subsection{Demonstration with superquadrics} \label{sec:demo:sq}

We now examine more general types of bodies, with superquadrics being typical examples. The Minkowski sums can be parameterized by their un-normalized gradients. Theorem \ref{thm:gradient} as well as the transformed version in both 2D and 3D cases are demonstrated. Some essential expressions to be used in the general formulas in Eq. \eqref{eq:mink_sum_transform_1} and \eqref{eq:mink_sum_transform_2} are shown here.

\subsubsection{Planar case}
The implicit function of a planar superquadric is
$$
\Psi(\xx) = \left( \frac{x_1}{a} \right)^{2 / \epsilon} + \left( \frac{x_2}{b} \right)^{2 / \epsilon} = 1 \,.
$$
Then, the gradient vector can be parameterized by the surface point as
$$ 
(\nabla \Psi)(\xx) = \begin{pmatrix}
\frac{2}{\epsilon a} \left(\frac{x_1}{a}\right)^{2/\epsilon - 1} \vspace{1ex} \\
\frac{2}{\epsilon b} \left(\frac{x_2}{b}\right)^{2/\epsilon - 1} \end{pmatrix}
\, \doteq \mm = \begin{pmatrix}
m_1 \vspace{1ex} \\
m_2 \end{pmatrix} \,. 
$$
By direct computations, the inverse function can be solved in terms of $\mm$ as
\begin{equation}
\tilde{\tilde{{\bf f}}}(\mm) \,=\,
\begin{pmatrix}
a \, \left(\frac{2}{\epsilon a}\right)^{\frac{-\epsilon}{2-\epsilon}}
m_{1}^{\frac{\epsilon}{2-\epsilon}} \vspace{1ex} \\
b \, \left(\frac{2}{\epsilon b}\right)^{\frac{-\epsilon}{2-\epsilon}}
m_{2}^{\frac{\epsilon}{2-\epsilon}} 
\end{pmatrix} \,. 
\label{eq:2d_sq_f_m}
\end{equation}

Theorem \ref{thm:gradient} is used here to compute the closed-form expression of Minkowski sums by first obtaining the relationships between $\mm_2$ and $\uu_2$ as
\begin{equation}
\mm_2 = {\bf g}_2(\uu_2) = \begin{pmatrix}
\frac{2}{a_2 \epsilon_2} u_{21}^{2-\epsilon_2} \vspace{1ex} \\
\frac{2}{b_2 \epsilon_2} u_{22}^{2-\epsilon_2}.
\end{pmatrix} \,,
\label{eq:2d_sq_m_u}
\end{equation}
where $u_{2j} = {\bf u}_2 \cdot {\bf e}_j$ and ${\bf e}_j$ is the $j$-th unit Cartesian coordinate basis, and the un-bolded $a_2, b_2$ and $\epsilon_2$ denote the shape parameters of $B_2$. Observing that ${\bf g}_2(k \uu_2) = k^{2-\epsilon_2} {\bf g}_2(\uu_2)$, the sufficient condition for this mapping (i.e. Lemma \ref{lemma:gradient:suff_cond}) is satisfied. So we are able to express $\| \mm_2 \|$ as a function of $\mm_1$ using Eq. \eqref{eq:Phi_m1_general}. The reverse function of ${\bf g}_2$ as a function of $\mm_1$ is then explicitly computed as
\begin{equation}
{\bf g}_2^{-1}(\mm_1) = \begin{pmatrix}
\left( \frac{a_2 \epsilon_2}{2} m_{11} \right)^{1/(2-\epsilon_2)} \vspace{1ex} \\
\left( \frac{b_2 \epsilon_2}{2} m_{12} \right)^{1/(2-\epsilon_2)}
\end{pmatrix} \,,
\label{eq:2d_sq_u_m}
\end{equation}
where $m_{1j} = \mm_1 \cdot {\bf e}_j$ denotes the $j$-th entry of the gradient vector of $B_1$. Finally, substituting Eqs. \eqref{eq:2d_sq_m_u} and \eqref{eq:2d_sq_u_m} into Eq. \eqref{eq:Phi_m1_general} gives the closed-form expression of $\|\mm_2\|$. 

Furthermore, using Eq. \eqref{eq:mink_sum_transform_2}, we demonstrate the closed-form Minkowski sums under different linear transformations for the two superellipses, including \footnote{Note that the right parenthesized superscript denotes the different types of transformation, and the right subscript denotes the $i$-th body.}:
\begin{itemize}
\item pure rotations, i.e. $M^{(1)}_i = rot(\alpha_i) = 
\left( \begin{array}{cc}
\cos \alpha_i & -\sin \alpha_i \\ \sin \alpha_i & \cos \alpha_i
\end{array} \right) $, where $\alpha$ is the rotational angle ;
\item pure shear, i.e. $M^{(2)}_i = shear(s_i) = 
\left( \begin{array}{cc}
1 & s_i \\ 0 & 1
\end{array} \right)$, where $s$ is the magnitude of the shear transformation ;
\item a combination of them, i.e. $M^{(3)}_i = rot(\alpha_i) \, shear(s_i)$ \,.
\end{itemize}

\begin{figure*}
\centering
\subfloat[Canonical forms.]{\includegraphics[scale=0.13]{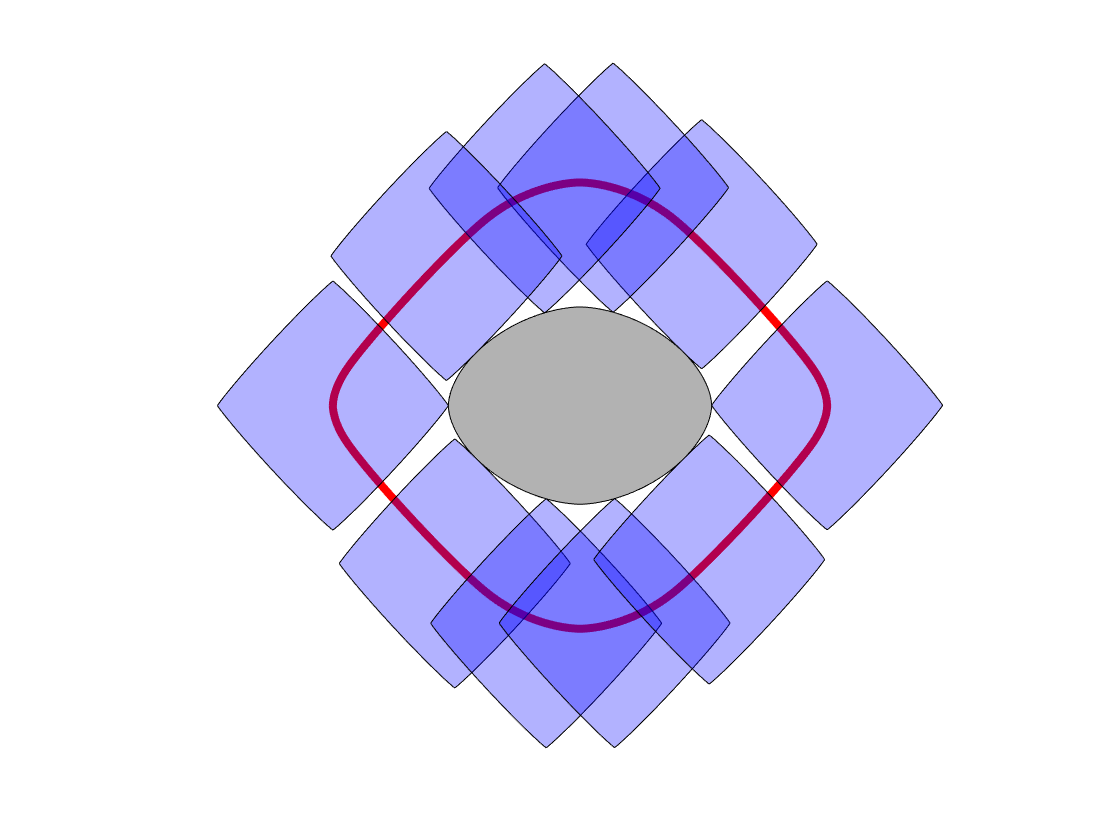}}
\subfloat[Pure rotations.]{\includegraphics[scale=0.13]{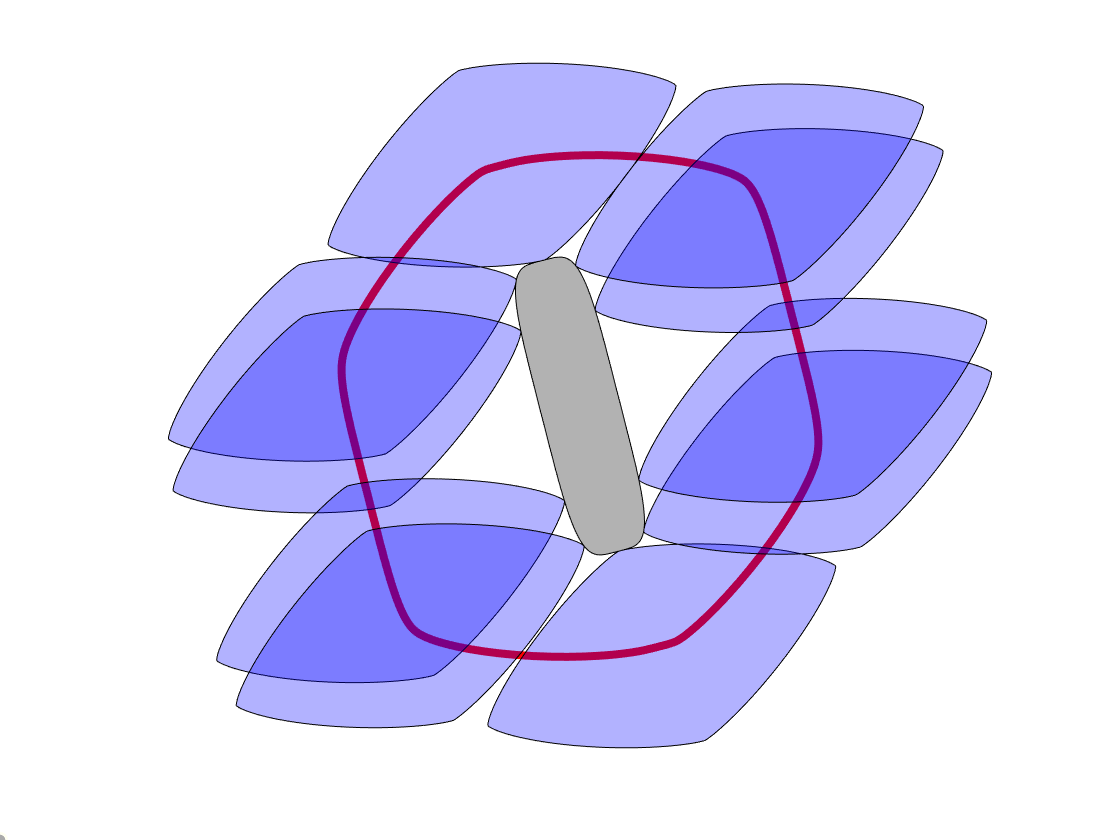}}
\subfloat[Pure shear transformations.]{\includegraphics[scale=0.13]{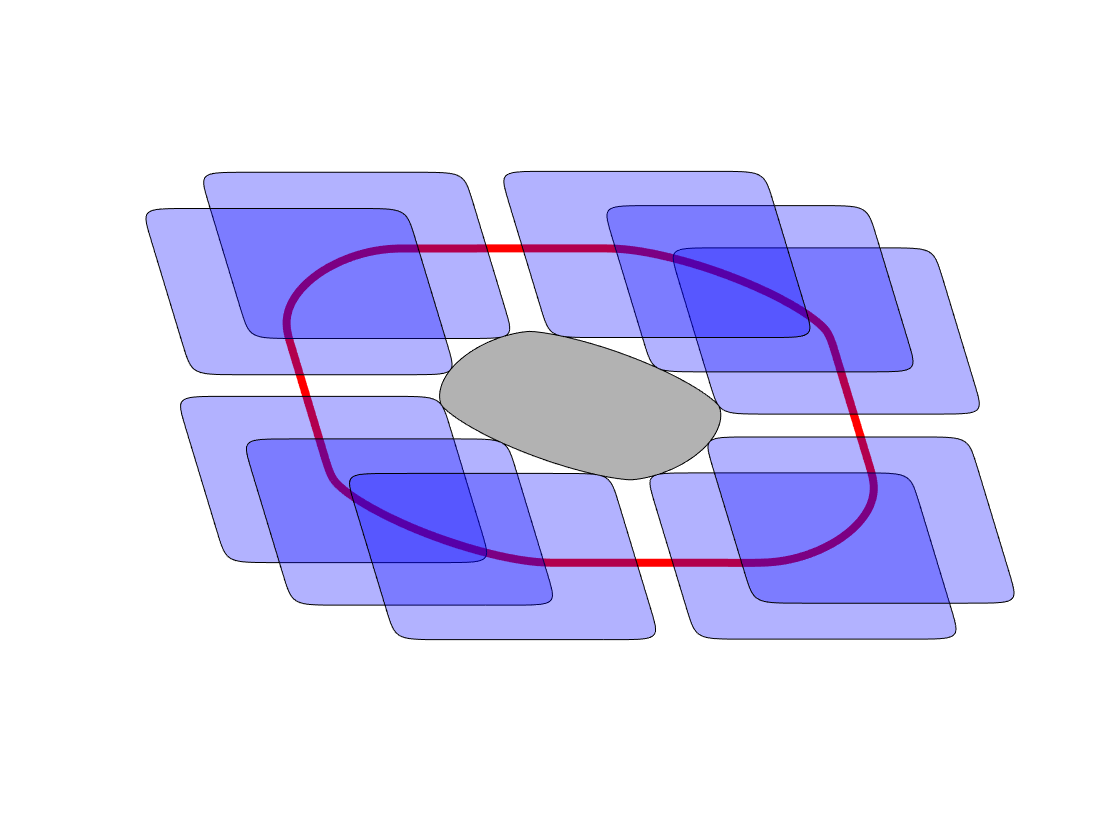}}
\subfloat[Combining rotation and shear.]{\includegraphics[scale=0.13]{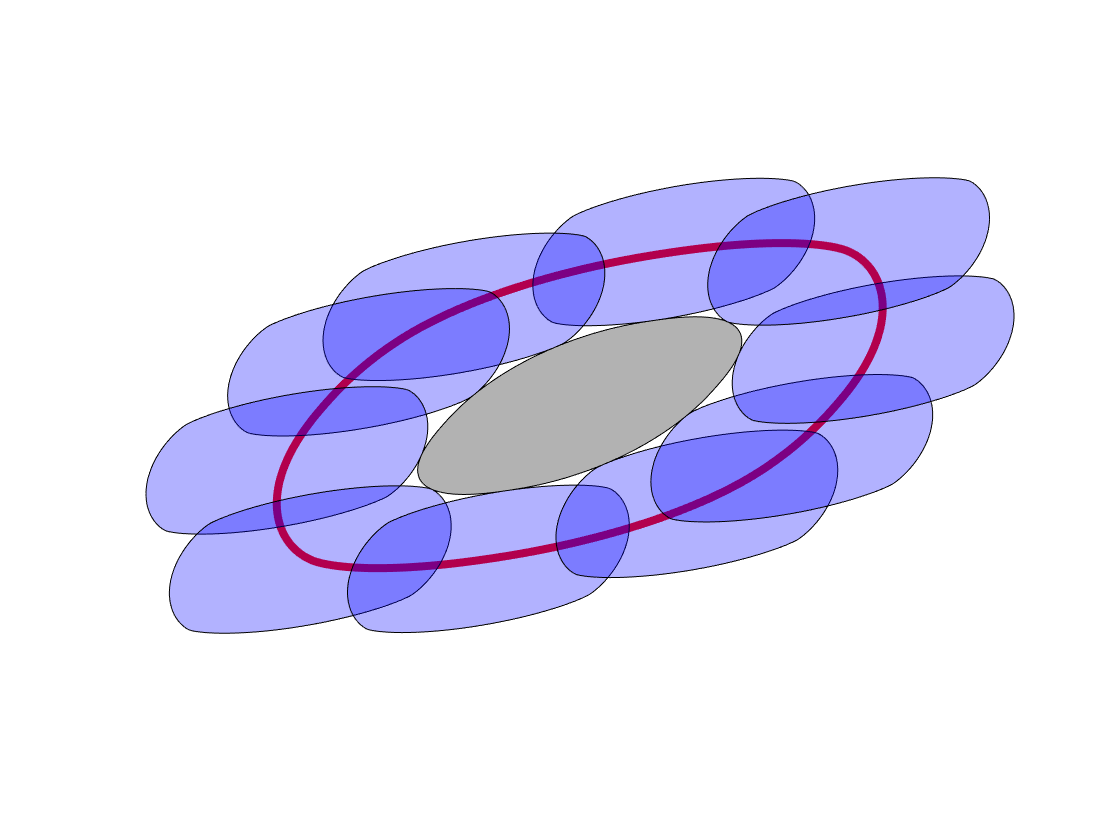}}
\caption{Demonstration of the closed-form Minkowski sums for two 2D superquadrics in canonical forms and under linear transformations. $B_1$ is plotted in grey color at the center. The red curve represents the computed closed-form Minkowski sums boundary. $B_2$ is plotted in blue color and placed at 10 different locations on the Minkowski sums boundary.}
\label{fig:2d_demo_sq}
\end{figure*}

Figure \ref{fig:2d_demo_sq} demonstrates the expressions using different types of convex superellipses with $\epsilon \in (0,2)$. $B_1$ is drawn in the center as a black curve, $B_2$ is drawn in blue and translates along the closed-form Minkowski sums boundary (in red) so as to kiss $B_1$. The demonstrations include superellipses in their canonical forms and under different types of linear transformations such as rotations and shear transformations.

\subsubsection{Spatial case}
The implicit expression for a 3D superquadric surface is
$$
\Psi(\xx) = \left( \left( \frac{x_1}{a} \right)^{2 / \epsilon_2} + \left( \frac{x_2}{b} \right)^{2 / \epsilon_2} \right)^{\epsilon_2 / \epsilon_1} + \left( \frac{x_3}{c} \right)^{2 / \epsilon_1} = 1,
$$
and the gradient is
$$
(\nabla \Psi) (\xx) =
\begin{pmatrix}
\frac{2}{a \epsilon_1} \left[ \psi(x_1, x_2) \right]^{\epsilon_2 / \epsilon_1 - 1} \left( \frac{x_1}{a} \right)^{2 / \epsilon_2 - 1} \vspace{1ex} \\
\frac{2}{b \epsilon_1} \left[ \psi(x_1, x_2) \right]^{\epsilon_2 / \epsilon_1 - 1} \left( \frac{x_2}{b} \right)^{2 / \epsilon_2 - 1} \vspace{1ex} \\
\frac{2}{c \epsilon_1} \left( \frac{x_3}{c} \right)^{2 / \epsilon_1 - 1} 
\end{pmatrix} \doteq \mm = 
\begin{pmatrix}
m_1 \\ m_2 \\ m_3
\end{pmatrix},
$$
where $\psi(x_1, x_2) = \left( \frac{x_1}{a} \right)^{2 / \epsilon_2} + \left( \frac{x_2}{b} \right)^{2 / \epsilon_2}$ denotes the part in the implicit expression that includes coordinates $x_1$ and $x_2$, and $\epsilon_1, \epsilon_2 \in (0,2)$ to ensure the convexity of the body enclosed by the bounding surface. By direct computations, the inverse function of gradient can be expressed as
\begin{equation}
\Tilde{ \Tilde{\ff} }(\mm) = 
\begin{pmatrix}
a \left( \frac{a \epsilon_1}{2} m_1 \right)^{\epsilon_2 / (2-\epsilon_2)} \left[ \gamma(m_3) \right]^{(\epsilon_1 - \epsilon_2) / (2-\epsilon_2)} \vspace{1ex} \\
b \left( \frac{b \epsilon_1}{2} m_2 \right)^{\epsilon_2 / (2-\epsilon_2)} \left[ \gamma(m_3) \right]^{(\epsilon_1 - \epsilon_2) / (2-\epsilon_2)} \vspace{1ex} \\
c \left( \frac{c \epsilon_1}{2} m_3 \right)^{\epsilon_1 / (2-\epsilon_1)}
\end{pmatrix} \,,
\label{eq:3d_sq_f_m}
\end{equation}
where 
$$
\gamma(m_3) = 1 - \left( \frac{c \epsilon_1}{2} m_3 \right)^{2 / (2-\epsilon_1)}.
$$

On the other hand, essential steps to compute $\Phi(\mm_1)$ includes
\begin{equation}
\mm_2 = {\bf g}_2(\uu_2) = 
\begin{pmatrix}
\frac{2}{a_2 \epsilon_{21}} \, u_{21}^{2-\epsilon_{22}} \left( u_{21}^2 + u_{22}^2 \right)^{(\epsilon_{22} - \epsilon_{21})/2} \vspace{1ex} \\
\frac{2}{b_2 \epsilon_{21}} \, u_{22}^{2-\epsilon_{22}} \left( u_{21}^2 + u_{22}^2 \right)^{(\epsilon_{22} - \epsilon_{21})/2} \vspace{1ex} \\
\frac{2}{c_2 \epsilon_{21}} \, u_{23}^{2-\epsilon_{21}}
\end{pmatrix} \,,
\label{eq:3d_sq_m_u}
\end{equation}
where $a_2, b_2, c_2, \epsilon_{21}$ and $\epsilon_{22}$ are the shape parameters of $B_2$, and $u_{2j} = \uu_2 \cdot {\bf e}_j$ is the $j$-th entry of the vector $\uu_2$ of $B_2$. It can also be observed that, ${\bf g}_2(k \, \uu_2) = k^{2-\epsilon_{21}} {\bf g}_2(\uu_2)$, which satisfies Lemma \ref{lemma:gradient:suff_cond}. Then the reverse function ${\bf g}^{-1}_2(\mm_1)$ can be computed as
\begin{equation}
{\bf g}_2^{-1}(\mm_1) = 
\begin{pmatrix}
\left( \frac{a_2 \epsilon_{21}}{2} m_{11} \right)^{\frac{1}{2-\epsilon_{22}}} \, \left[ \rho(m_{21}, m_{22}) \right]^{\frac{\epsilon_{21}-\epsilon_{22}}{4-2 \epsilon_{21}}} \\
\left( \frac{b_2 \epsilon_{21}}{2} m_{12} \right)^{\frac{1}{2-\epsilon_{22}}} \, \left[ \rho(m_{21}, m_{22}) \right]^{\frac{\epsilon_{21}-\epsilon_{22}}{4-2 \epsilon_{21}}} \\
\left( \frac{c_2 \epsilon_{21}}{2} m_{13} \right)^{\frac{1}{2-\epsilon_{21}}}
\end{pmatrix} \,,
\label{eq:3d_sq_u_m}
\end{equation}
where 
$$
\rho(m_{21}, m_{22}) = \left( \frac{a_2 \epsilon_{21}}{2} m_{21} \right)^{\frac{2}{2-\epsilon_{22}}} + \left( \frac{b_2 \epsilon_{21}}{2} m_{22} \right)^{\frac{2}{2-\epsilon_{22}}} \,,
$$ 
and $m_{ij} = \mm_i \cdot {\bf e}_j$ is the $j$-th entry of the gradient vector $\mm_i$ of $B_i$.

\begin{figure*}
\centering
\subfloat[Two general superquadrics in canonical forms. \label{fig:3d_demo_sq:sq_sq_canonical}]{\includegraphics[scale=0.3, trim={50 50 50 40},clip]{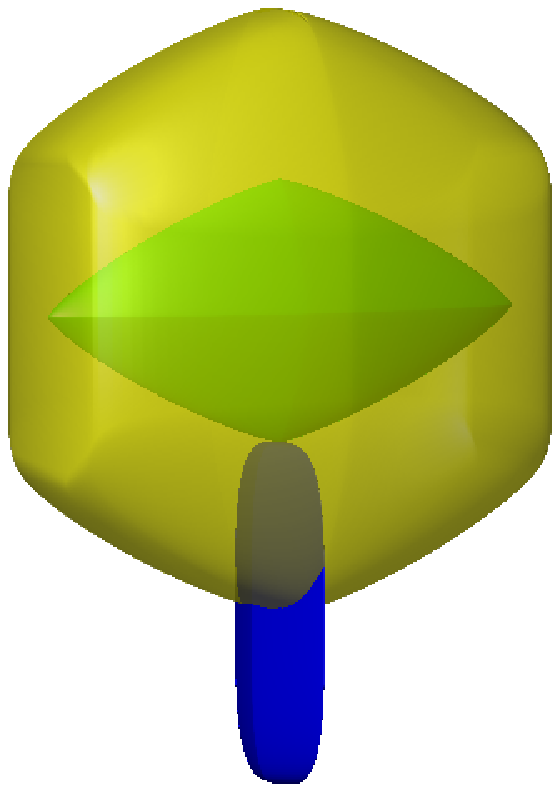}} 
~
\subfloat[Two general superquadrics with pure rotations. \label{fig:3d_demo_sq:sq_sq_rot}]{\includegraphics[scale=0.3, trim={50 50 50 50},clip]{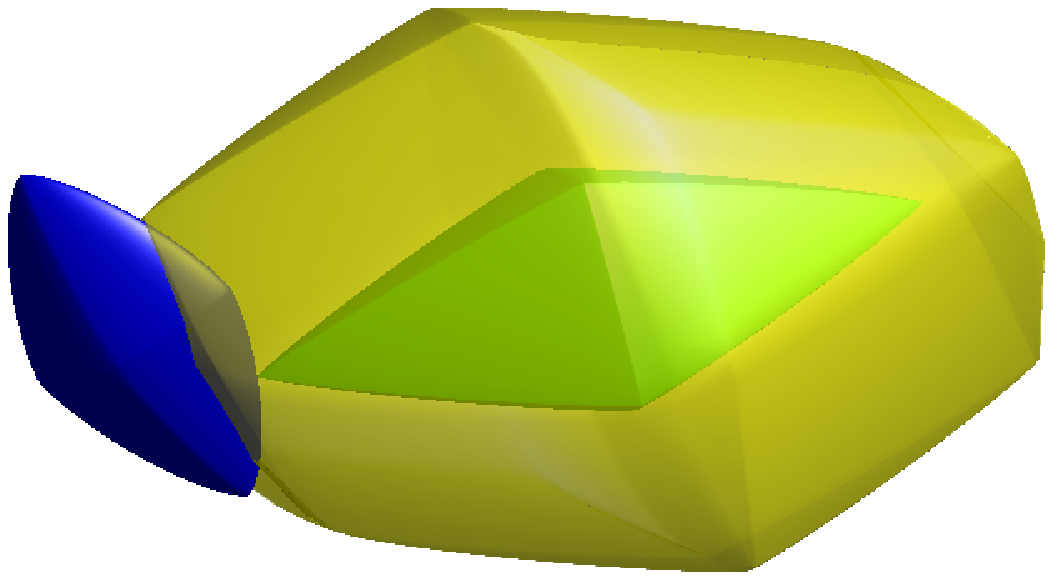}}
~
\subfloat[Two cubes. \label{fig:3d_demo_sq:cube_cube}]{\includegraphics[scale=0.3, trim={50 50 60 30},clip]{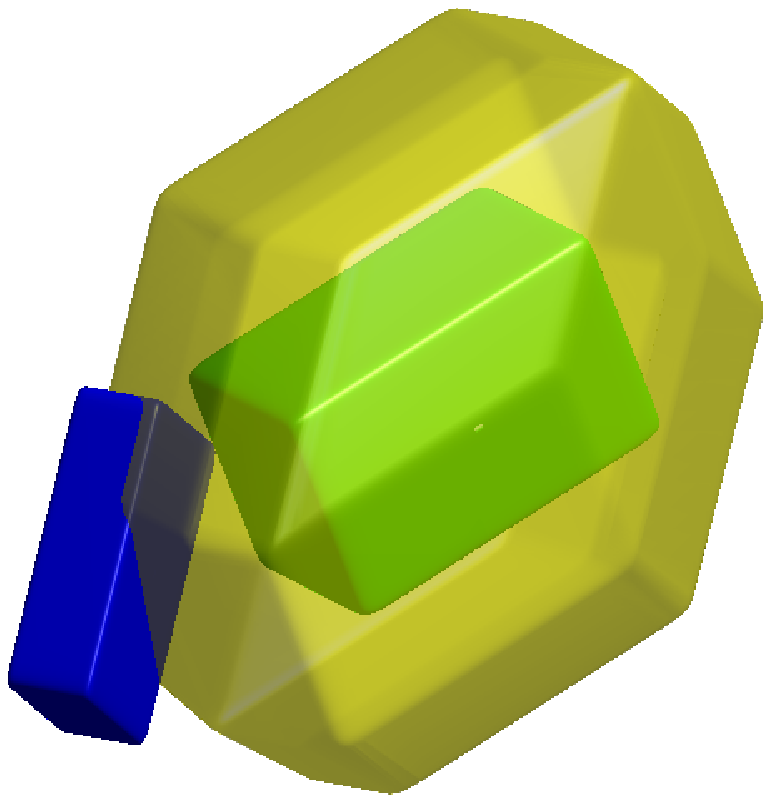}}
~
\subfloat[Cube and elliptical cylinder. \label{fig:3d_demo_sq:cube_cyl}]{\includegraphics[scale=0.3, trim={50 50 50 50},clip]{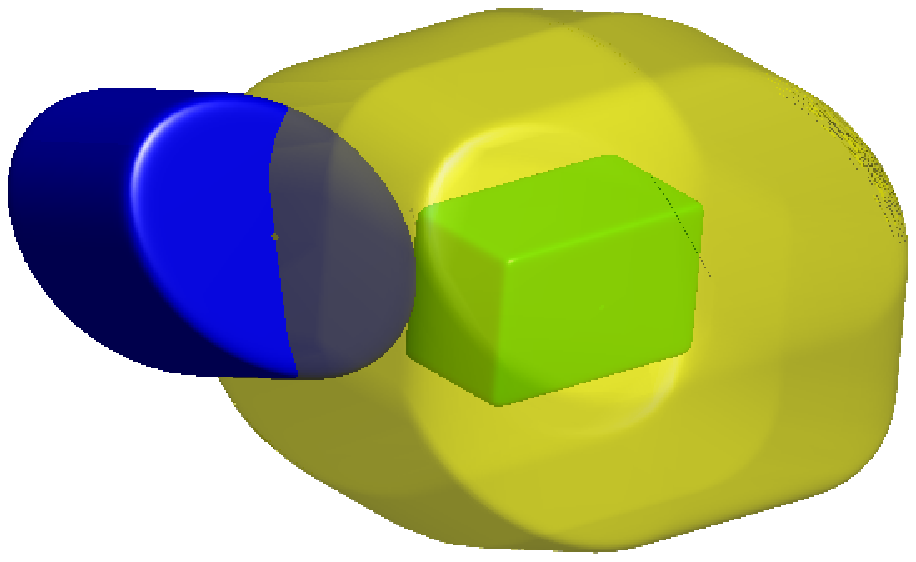}}
\\
\subfloat[Two parallelepipeds. \label{fig:3d_demo_sq:parpipde_parpiped}]{\includegraphics[scale=0.3, trim={50 50 50 50},clip]{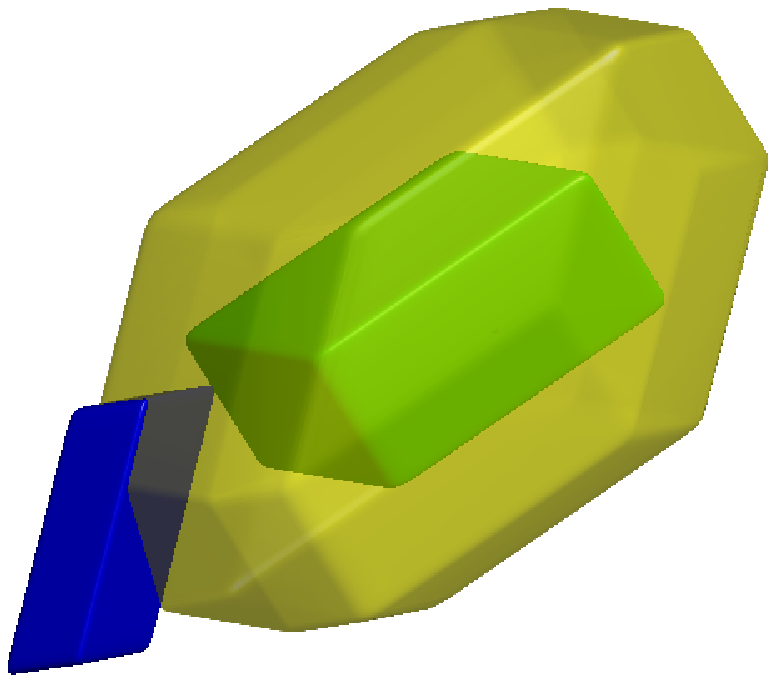}}
~
\subfloat[Parallelepiped and elliptical cylinder. \label{fig:3d_demo_sq:parpiped_cyl}]{\includegraphics[scale=0.3, trim={60 80 60 80},clip]{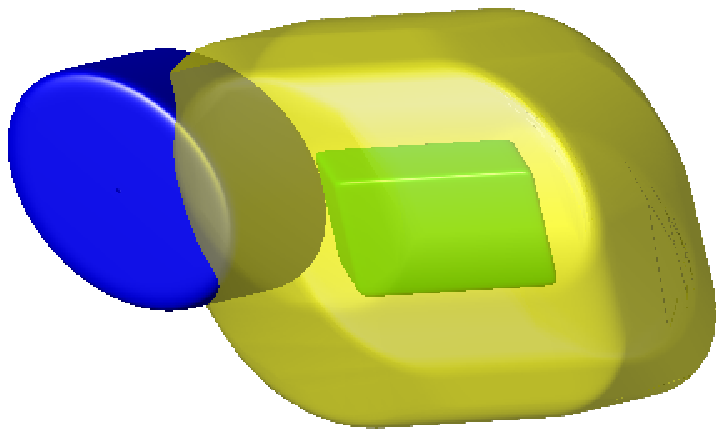}}
~
\subfloat[Two elliptical cylinders. \label{fig:3d_demo_sq:cyl_cyl}]{\includegraphics[scale=0.3, trim={50 50 50 50},clip]{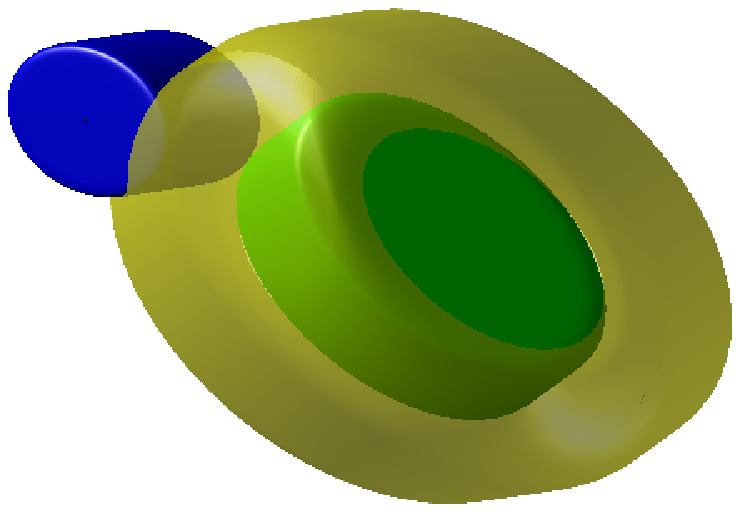}}
~
\subfloat[Two ellipsoids. \label{fig:3d_demo_sq:ellip_ellip}]{\includegraphics[scale=0.3, trim={50 50 50 50},clip]{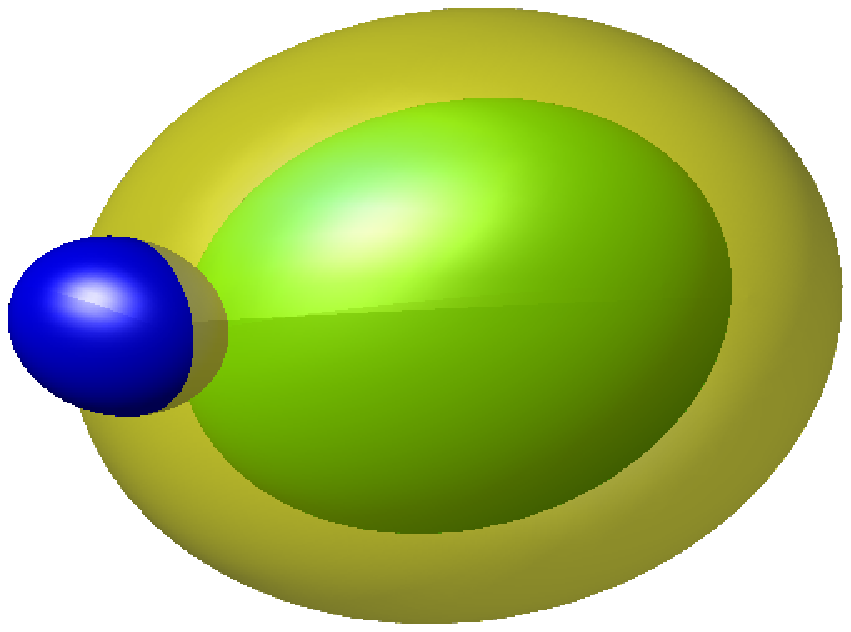}}
\caption{Demonstrations of the closed-form Minkowski sums for two 3D superquadric models in general forms and representing basic geometric primitives. $B_1$ is shown in green color at the center. The yellow surface represents the computed closed-form Minkowski sums boundary. $B_2$ is plotted in blue color and placed at one location on the Minkowski sums boundary.}
\label{fig:3d_demo_sq}
\end{figure*}

The closed-form Minkowski sums for $B_1$ and $B_2$ in canonical form and with different orientations are visualized here, as shown in Figs. \ref{fig:3d_demo_sq:sq_sq_canonical} and \ref{fig:3d_demo_sq:sq_sq_rot} repectively. $B_1$ is drawn in the center in green, and $B_2$ is drawn in blue and translates along the closed-form Minkowski sums boundary, which only touches $B_1$ at every point.

\subsubsection{Examples on basic geometric primitives} \label{sec:demo:sq:geom_primitives}
Basic geometric primitives, like cube, parallelepiped and cylinder, are fundamentals on computer-aided design, and superquadrics can either exactly represent or approximate them well. Here, concrete examples of Minkowski sums for these primitives by varying superquadric parameters are demonstrated. In the following contents, the semi-axes lengths and the exponents of superquadrics are denoted as a tuple of $(a, b, c)$ and a pair of $(\epsilon_1, \epsilon_2)$, respectively.

\begin{itemize}
\item Cube. A cube can be defined by a tuple of length, width and height, i.e., $(l, w, h)$. It can be approximated by a superquadric model as: $(a, b, c) = (l/2, w/2, h/2)$ and $(\epsilon_1, \epsilon_2) = (0.1, 0.1)$.

\item Parallelepiped. A parallelepiped can be seen as a linearly transformed cube. The linear transformation is defined in the matrix $M$. Then the superquadric approximation of the parallelepiped uses the same parameters with the cube and transformed by $M$.

\item Ellipsoid. An ellipsoid is defined by semi-axes length $(a_{e}, b_{e}, c_{e})$ and can be defined exactly using the superquadric model: $(a, b, c) = (a_{e}, b_{e}, c_{e})$ and $(\epsilon_1, \epsilon_2) = (1, 1)$. When $a_{e} =  b_{e} = c_{e} = r$, it becomes a sphere with radius $r$. Note that, in the case of two ellipsoids, their closed-form Minkowski sums are expressed in simpler forms of Eqs. \eqref{eq:ellipsoid_mink_sum_thm_1} and \eqref{eq:ellipsoid_mink_sum_closed_m}.

\item Cylinder. A general elliptic cylinder is defined by the elliptical cross section with semi-axes $(a_{c}, b_{c})$ and height $h$. The approximated superquadric parameters are then defined as: $(a, b, c) = (a_{c}, b_{c}, h/2)$ and $(\epsilon_1, \epsilon_2) = (0.1, 1)$.
\end{itemize}

Using the above superquadric representations, Figs \ref{fig:3d_demo_sq:cube_cube}, \ref{fig:3d_demo_sq:cube_cyl}, \ref{fig:3d_demo_sq:parpipde_parpiped}, \ref{fig:3d_demo_sq:parpiped_cyl}, \ref{fig:3d_demo_sq:cyl_cyl} and \ref{fig:3d_demo_sq:ellip_ellip} demonstrate the Minkowski sums using the proposed closed-form expression of Eq. \eqref{eq:mink_sum_transform_2} for different pairs of basic geometric primitives. 

\section{Numerical Validations and Benchmarks in Simulation for the Closed-form Minkowski Sums} \label{sec:numerial_verifications}
To validate the correctness and evaluate the performance of the proposed closed-form Minkowski sums, we conduct numerical simulations and compare with the discrete Minkowski sums computations from definition as shown in Fig. \ref{fig:mink:demo_mink_sum:def}. 

At first, a set of points on the Minkowski sums boundary is generated using the proposed exact closed-form expression. And we examine, when placing $B_2$ at these points, whether there exists a kissing point that locates on both body boundary surfaces and the corresponding normal vectors are anti-parallel. To evaluate the efficiency of the proposed closed-form expression, we conduct benchmarks in the second simulation. The numbers of points on the boundaries of $S_1$ and $S_2$ are varied, and the running time to generate the point set on the Minkowski sums boundary are compared. All the simulations throughout this article are implemented in Matlab R2020b and performed in an Intel Xeon CPU at 2.80 GHz.

\subsection{Numerical validations of the closed-form Minkowski sums computations}

The general idea is to compute the kissing point of the two bodies when placing $B_2$ at each point on Minkowski sum boundary. We are using the fact that the kissing point $\xx_{kiss}$ between two bodies must satisfies: 

\begin{itemize}
\item The kissing point lies on both surfaces. This can be achieved when the implicit expression at the kissing point equals to one, i.e. 
\begin{equation}
\Psi_1 \left( ^{1}\xx_{kiss} \right) = \Psi_2 \left( ^{2}\xx_{kiss} \right) = 1,
\label{eq:kiss_point}
\end{equation}
where $\Psi_i(\cdot)$ is the implicit expression of the $i$-th body and $^{i}\xx_{kiss}$ is the kissing point as viewed in the body frame of $i$-th body.

\item The normal vectors at the kissing point with respect to the two surfaces are anti-parallel. Specifically, we can compute the un-normalized gradients at the kissing point on each surface, and verify that the angle between them is $\pi$.
\end{itemize}

In practice, we first sample a set of points $\xx_1 \left(\pphi_k \right)$ on $\partial B_1$ based on the set of angular parameters $\left\{ \pphi_k \right\}$. Then, this parameter set can also be used to generate the points on Minkowski sum $\xx_{1+\overline{2}}(\pphi_k)$ using our closed-form expression. Then, we place $B_2$ at $\xx_{1+\overline{2}} \left( \pphi_k \right)$, and verify that each $\xx_1 \left( \pphi_k \right)$ is the kissing point. 

To validate the first condition that the kissing point locates on both bounding surfaces, $\xx_1(\pphi_k)$ is transformed into the local frame of each body, then Eq. \eqref{eq:kiss_point} can be examined. Note that, since $\xx_1$ is computed from the parametric expression of $B_1$, we only need to verify that the value of
\begin{equation}
e_{implicit} \doteq \left| \Psi_2(^{2} \xx_1) - 1 \right| = 0.
\label{eq:verify_implicit}
\end{equation}
Then, since the gradient can be computed from the point on surface, the un-normalized gradients at the candidate kissing point on both bodies can be computed as $\mm_1 \left( \xx_1 \right)$ and $\mm_2 \left( \xx_1 \right)$ respectively. They are anti-parallel when
\begin{equation}
e_{gradient} \doteq \left| \frac{\mm_1(\xx_1) \cdot \mm_2(\xx_1)}{\|\mm_1(\xx_1)\| \|\mm_2(\xx_1)\|} + 1 \right| = 0.
\label{eq:verify_gradient}
\end{equation}
Once these two conditions are verified, we can say that $\xx_1$ is the kissing point when $B_2$ is place at the Minkowski sum boundary.

We conduct simulations to verify in both 2D and 3D cases. For 2D case, a total of 1000 points are sampled on $\partial B_1$ based on the angular parameter set $\left\{ \pphi^{2D}_k \right\} = \left\{ \theta_k \right\}$; and for 3D case, 100 samples are generated for each angular parameter pair, i.e. $\left\{ \pphi^{3D}_k \right\} = \left\{ (\eta_k, \omega_k) \right\}$, so a total of $10^{4}$ points are sampled on $\partial B_1$. And for both cases, 100 trials are simulated and the mean values of Eq. \eqref{eq:verify_implicit} and \eqref{eq:verify_gradient} for among all the trials are computed.

For the 2D case, among all the trials, the mean values of $e_{implicit}$ and $e_{gradient}$ are $5.9194 \times 10^{-16}$ and $6.1839 \times 10^{-7}$, respectively. And for the 3D case, the mean values of $e_{implicit}$ and $e_{gradient}$ are $1.8303 \times 10^{-9}$ and $4.0654 \times 10^{-7}$, respectively. The results show that when $B_2$ is centered at the Minkowski sum boundary points $\xx_{1+\overline{2}}(\pphi)$, which are computed by both our proposed closed-form expression, there exists a kissing point $\xx_1 \left( \pphi \right)$. Ideally, the expressions in Eq. \eqref{eq:verify_implicit} and \eqref{eq:verify_gradient} should equal to zero since our proposed closed-form Minkowski sums expression is the exact solution. Due to the numerical precision issues, there are discrepancies between the computational results and the ideal expectations. However, these discrepancies are all very small, i.e. less than the level of $10^{-6}$, which is acceptable for the numerical verification purposes. Therefore, the point on $\partial B_1$ parameterized by $\pphi_k$ is the kissing point when placing $B_2$ at $\xx_{1+\overline{2}} \left( \pphi_k \right)$, with which the proposed closed-form Minkowski sum expression is numerical verified.

\subsection{Fitting approximation errors of different geometric models}

The Minkowski sums expressions proposed here are \emph{exact} provided that the bodies are convex and enclosed by smooth positively curved boundaries. The superquadric model studied in Sec. \ref{sec:demo:sq} is one of the typical examples, which is flexible in varying shapes with only a small amount of parameters. Alternatively, one may consider the body to be bounded by discrete surfaces, e.g., a convex polyhedra. A convex polyhedra is a well-known type of geometry to represent a complex body, but may require many vertices and faces to approximate a rounded region. It is possible to fit a convex polyhedra with superquadric model and vice versa, which introduces approximation errors. Therefore, to evaluate the approximation quality, the fitting errors of different types of geometric models are studied in this subsection. The evaluation metric is based on the relative volume error of the approximated body with respect to the ground truth, i.e.,
\begin{equation}
e_{volume} = \frac{\left| V(B_{query}) - V(B_{true}) \right|}{V(B_{true})} \,,
\label{eq:fit_error_volume}
\end{equation}
where $V(\cdot)$ denotes the volume of a body, $B_{query}$ and $B_{true}$ denotes the queried body and reference body which is treated as ground truth, respectively. In the following studies, for each case, a total number of 100 trials are simulated.

\subsubsection{Case 1: fitting superquadrics to convex polyhedra}
For a given convex polyhedron in general position with $m$ vertices, a superquadric fitting model can be computed by minimizing \citep{vaskevicius2017revisiting}
\begin{equation}
\min_{a,b,c,\epsilon_1,\epsilon_2,R,{\bf t}} abc \sum_{k=1}^{m} \left[ \Psi(R^T(\xx_{k}-{\bf t}))^{\epsilon_1} - 1 \right]^2 \,,
\label{eq:fit_sq_to_points}
\end{equation}
where $\xx_{k}$ are the vertices of the polyhedron, $(a,b,c)$, $(\epsilon_1, \epsilon_2)$, $R$ and ${\bf t}$ are optimization variables denoting semi axes length, exponent, orientation and center position of the superquadric model, respectively. The optimization attempts to approach the superquadric surface to all the vertices of the convex polyhedron and minimize the volume of the region it encloses using the term $abc$. The volume of a superquadric body can be computed as \citep{jaklivc2000superquadrics} 
$$ 
V(B_{sq}) = 2 \, abc \, \epsilon_1 \epsilon_2 \, \beta \left( \frac{\epsilon_1}{2} + 1, \epsilon_1 \right) \, \beta \left( \frac{\epsilon_2}{2}, \frac{\epsilon_2}{2} \right) \,,
$$ 
where $\beta(x,y) = 2 \int_0^{\pi/2} \sin^{2x-1} \phi \, \cos^{2y-1} \phi \, d \phi$ is the beta function.

In this case, the convex polyhedron is treated as ground truth and generated as the convex hull of randomly point set. The approximation error metrics include not only Eq. \eqref{eq:fit_error_volume}, but also the absolute difference between the point and implicit function, i.e., Eq. \eqref{eq:verify_implicit}. In each simulation trial, a set of 100 random points are sampled, whose convex hull has an average of 30 vertices. Among all the trials, the mean of $e_{volume}$ and $e_{implicit}$ are $0.2095$ and $0.5808$, respectively. The results show that the superquadric surfaces fit closely to the polyhedral vertices. But it should also note that when the convex polyhedron is highly non-symmetric, fitting the central-symmetric superquadric model is conservative and the volume difference might be non-avoidably large. More detailed analysis of a superquadric fitting problem is out of the scope of this article, but can be referred to \citep{vaskevicius2017revisiting}.

\subsubsection{Case 2: fitting convex polyhedra to superquadrics} \label{sec:numerical_verifications:poly_sq}
This case considers the ground truth to be a superquadric body. The approximated convex polyhedron is defined by sampling points on the smooth superquadric surface. The sampling density significantly affects the approximation accuracy. Therefore, in this case, the relative volume is measured in different simulation trials by varying the number of sampled points. Among the 100 trials for each number, ranging from 16 to 400, the averaged $e_{volume}$ drops monotonically from $0.6573$ (16 sampled points) to less than $0.1$ (i.e., $0.094$ when 81 points are sampled), and finally reaches $0.0171$ (when 400 points are sampled). It can be seen that a good convex polyhedral approximation uses many more sampled points. This is mainly because that a better faceted representation of the curved surface of a superquadric requires denser set of sampled points.

\subsubsection{Case 3: fitting superquadrics and convex polyhedra to basic geometric primitives}
Both superquadrics and convex polyhedra are able to exactly or approximately represent the basic geometric primitives like cubes, ellipsoids and cylinders. This case compares the fitting errors for both of them. The superquadric models for these basic geometric primitives are introduced in Sec. \ref{sec:demo:sq:geom_primitives}.

\begin{itemize}
\item Cube. Convex polyhedron can exactly represent a cube by defining vertices as the corner points of the cube. The fitting error for superquadric model among the 100 trials is $0.0112$ on average.

\item Ellipsoid. Superquadric is an exact model for ellipsoid primitive, so there is no fitting error. The convex polyhedral body, on the other hand, approximates an ellipsoid by sampling points on the surface. In the simulation trials, 100 points are sampled. And the averaged fitting error is $0.1340$.

\item Cylinder. Both superquadric and convex polyhedron are approximations for a general elliptic cylinder. The averaged fitting error for a superquadric model is $0.0074$. Here, the number of vertices of a convex polyhedron is tuned to be $60$. With this number, the averaged fitting error achieves $0.0078$ among all the simulation trials, which is at the same level with a superquadric model.
\end{itemize}

The superquadric model approximates basic primitives which are central-symmetric, i.e., cube and cylinder, well. The relative volume errors with the true geometries are around or less that $1\%$. On the other hand, for the ellipsoidal model, even 100 points are used to define a convex polyhedron, the volume error is still as large as $10\%$. For the cylinder case, both superquadric and convex polyhedron are approximations. To achieve the same level of accuracy in volume, a convex polyhedron needs as many as $60$ vertices, as compared to $11$ parameters for a superquadric model.

\subsection{Benchmarks for computational time of the point sets on Minkowski sums boundary}

The computational time for generating the Minkowski sums boundary are measured for our proposed method and some other existing methods. Three other methods are compared: \emph{convex hull} \citep{lozano1990spatial}, \emph{edge sorting} \citep{mark2008computational} and \emph{closed-form geometric} \citep{yan2015closed, ruan2018path} methods. Three cases of body pairs are studied: ellipsoid-ellipsoid (\emph{EE}), ellipsoid-superquadrics (\emph{ES}) and superquadrics-superquadrics (\emph{SQ}). And both 2D and 3D spaces are considered and compared.

The \emph{convex hull} method is introduced in Sec. \ref{sec:general_mink:def} and only works for two convex bodies. The procedure is to firstly sample points on the smooth surface of the two bodies. Then, the summation of all the pairs of points from each body surface is computed, which results in a superset of the Minkowski sums. To finally obtain the boundary, convex hull is then computed. It works for any dimension as long as the convex hull computation is available. In the benchmark, we use the QuickHull \citep{barber1996quickhull} algorithm to compute the convex hull, which is the built-in implementation of convex hull in Matlab. The \emph{edge sorting} method achieves $O(m+n)$ complexity for two planar polygons. This method merges the slope of the edges of each body and add the vertices in a sorted order, making it run in linear time. Although both \emph{convex hull} and \emph{edge sorting} methods do not provide exact results in the case of smooth bodies, we could still get a sense of the speed when generating points on Minkowski sums boundary. For two bodies with smooth surfaces, the \emph{closed-form geometric} method solves for exact expressions of Minkowski sums also in closed-form. It works when at least one body is an ellipsoid. The method proposed in this article is denoted as \emph{closed-form (proposed)}. We use the more general gradient-based expression, i.e. Eq. \eqref{eq:mink_sum_transform_2}, for all benchmark simulations. Both closed-form methods work in any dimension.

\begin{figure*}
\centering
\subfloat[2D case]{\includegraphics[scale=0.12]{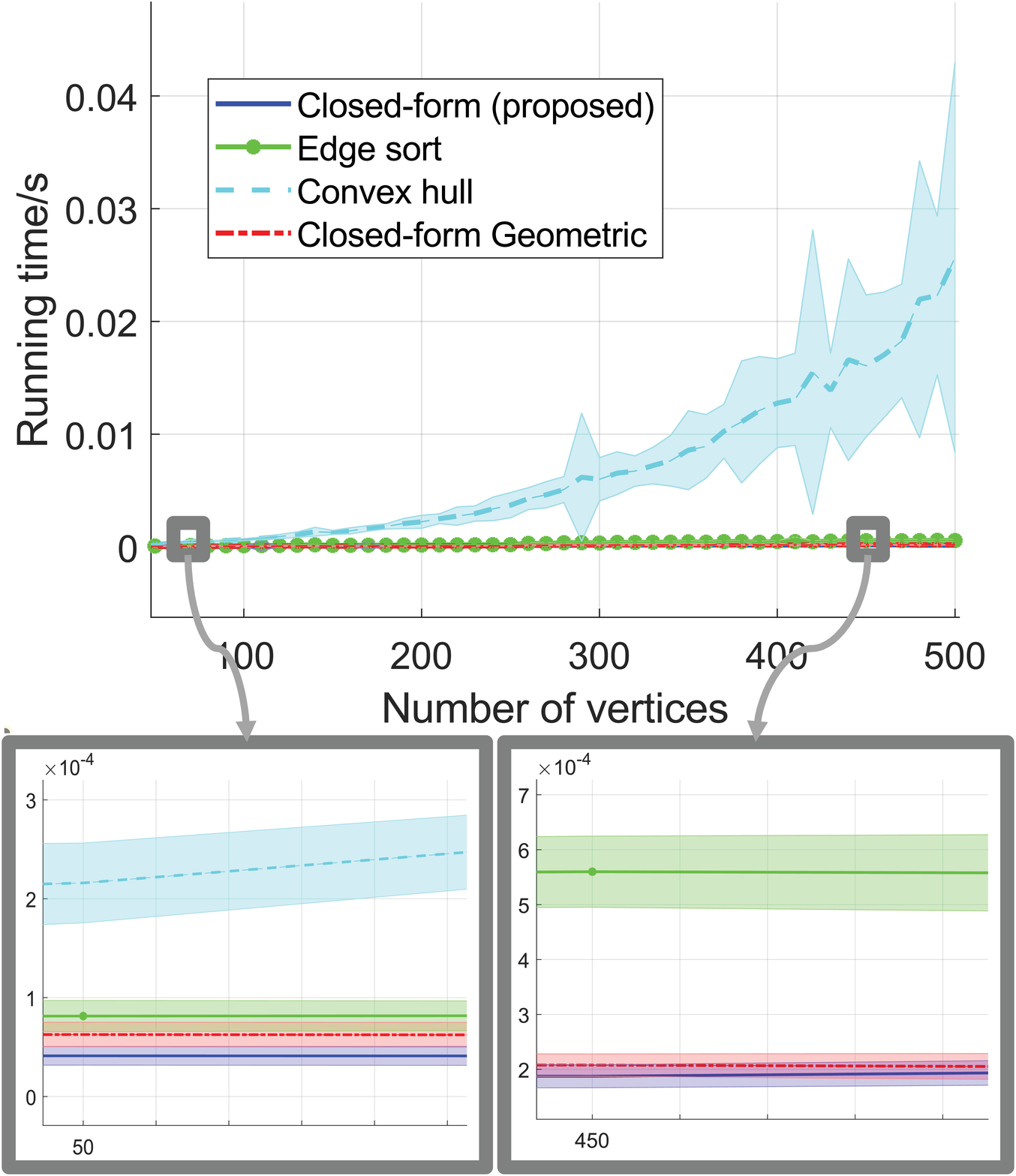}}
~
\subfloat[3D case]{\includegraphics[scale=0.12]{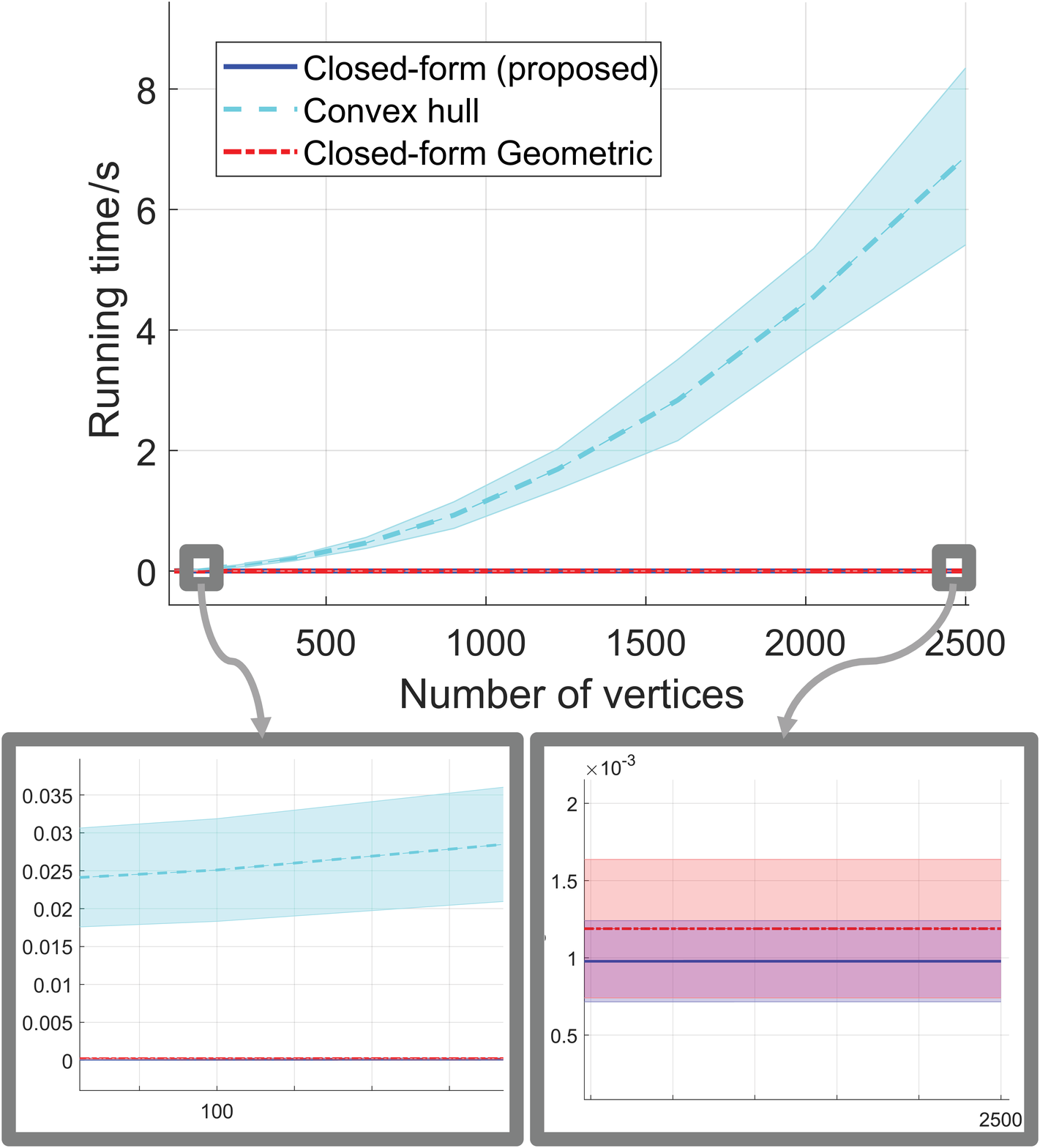}}
\caption{Running time comparisons for different algorithms in computing Minkowski sums between two ellipsoids (EE case). For all sub-figures, curves represents the mean value of the running time, and the shaded regions show the range of the standard deviation. Two zoom-in views of small regions are shown below each sub-figure.}
\label{fig:verification:run_time:ee}
\end{figure*}

\begin{figure*}
\centering
\subfloat[2D case]{\includegraphics[scale=0.12]{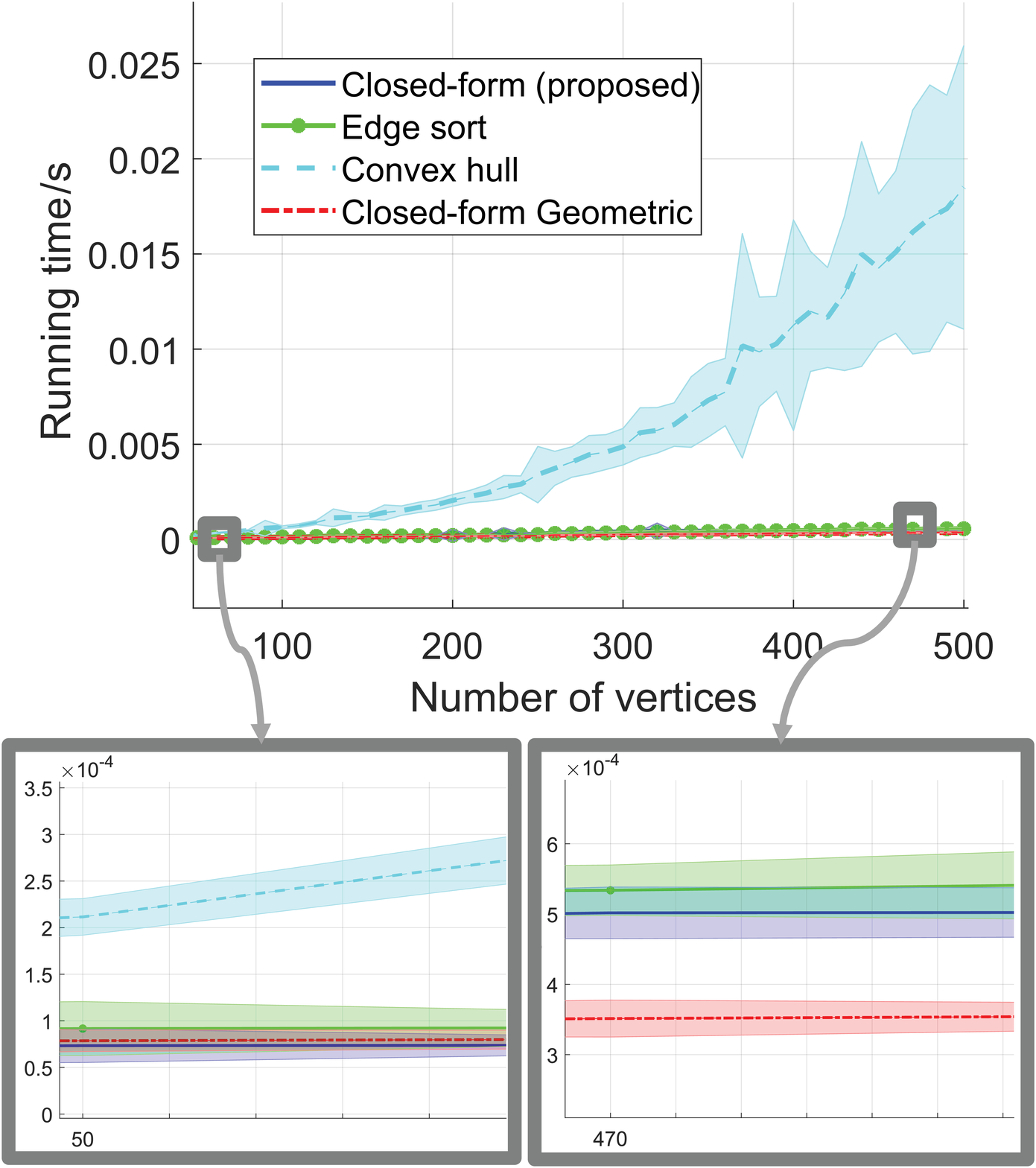}}
~
\subfloat[3D case]{\includegraphics[scale=0.12]{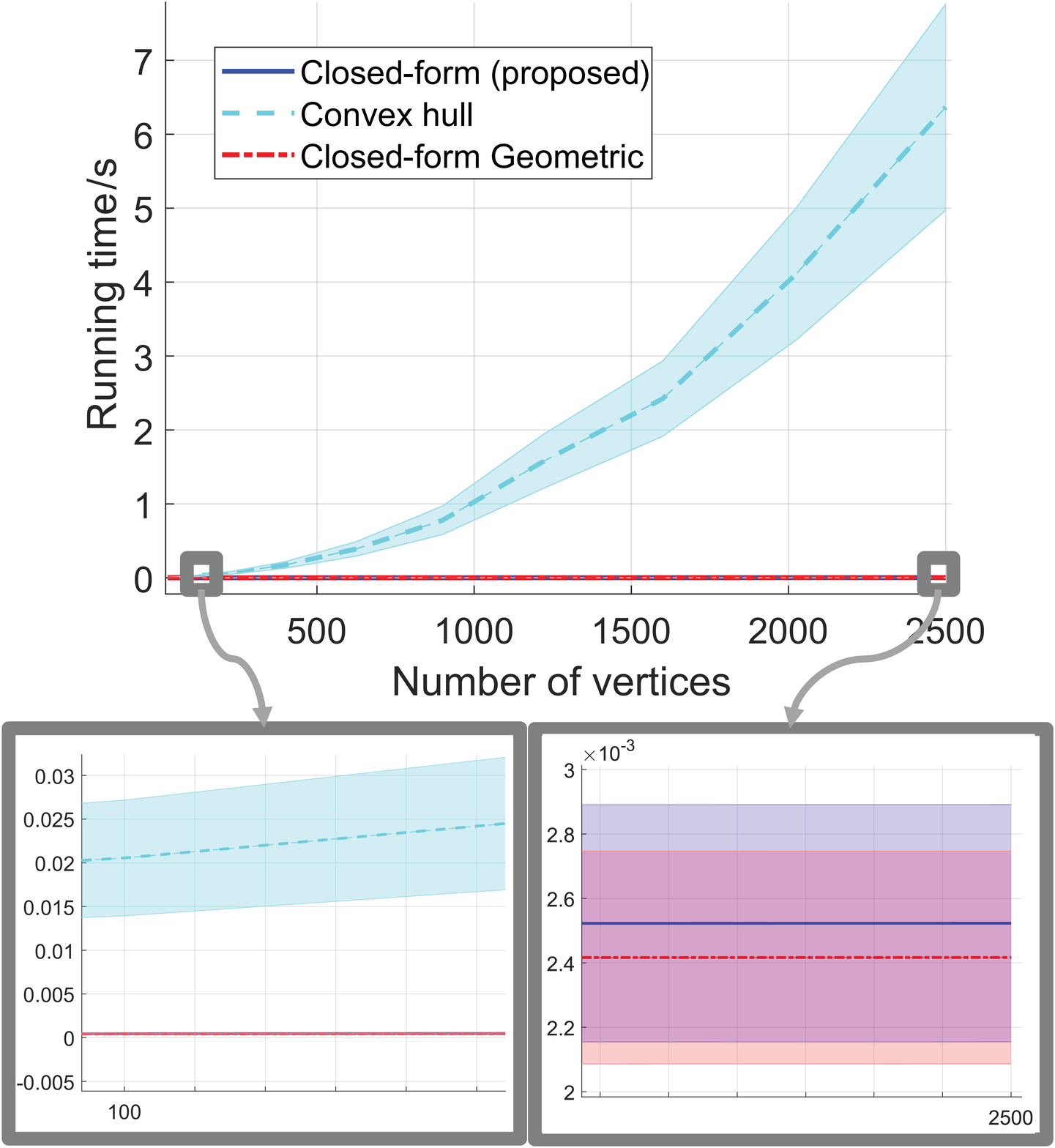}}
\caption{Running time comparisons for different algorithms in computing Minkowski sums between one ellipsoid and a superquadric (ES case).}
\label{fig:verification:run_time:es}
\end{figure*}

\begin{figure*}
\centering
\subfloat[2D case]{\includegraphics[scale=0.12]{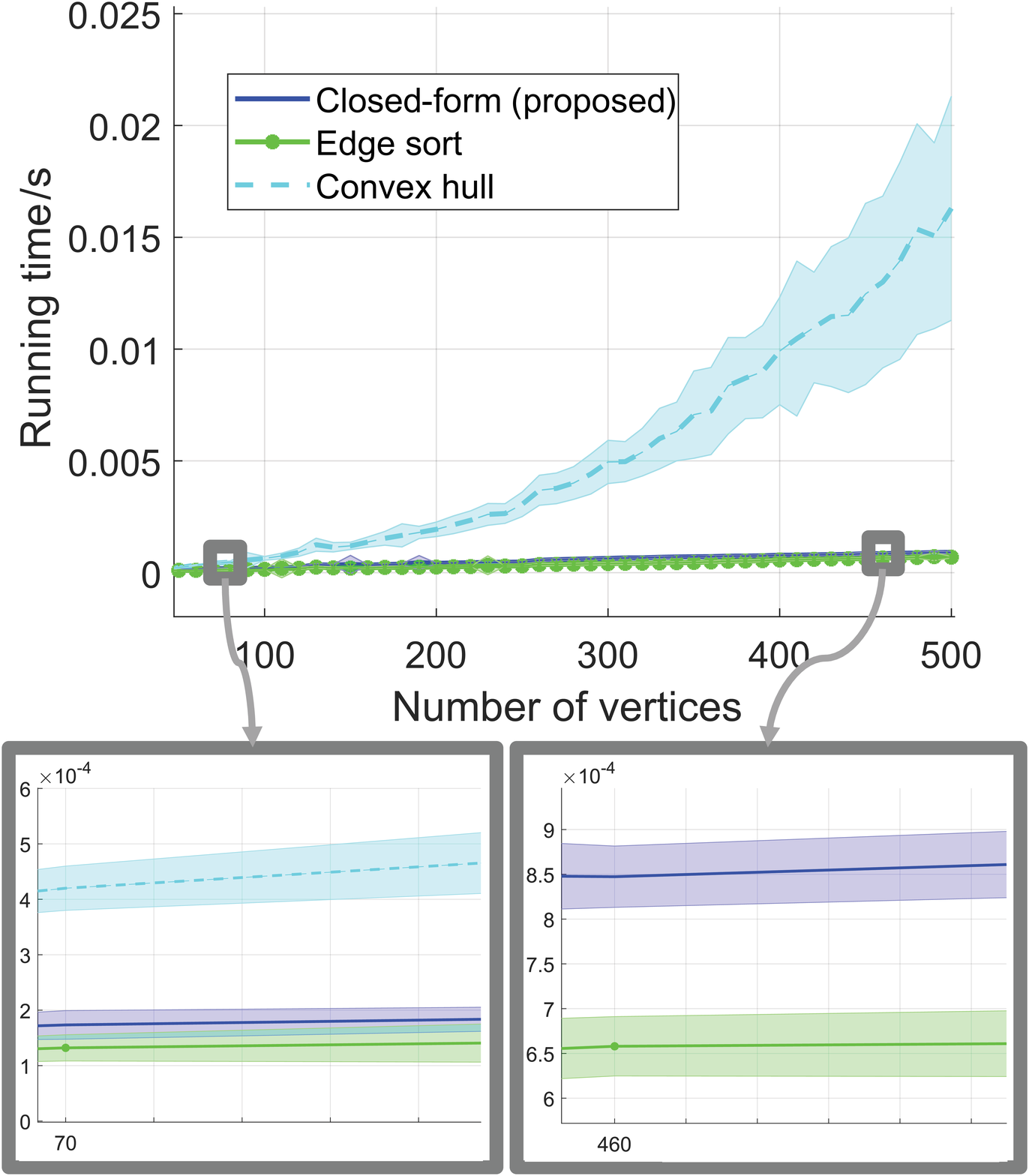}}
~
\subfloat[3D case]{\includegraphics[scale=0.12]{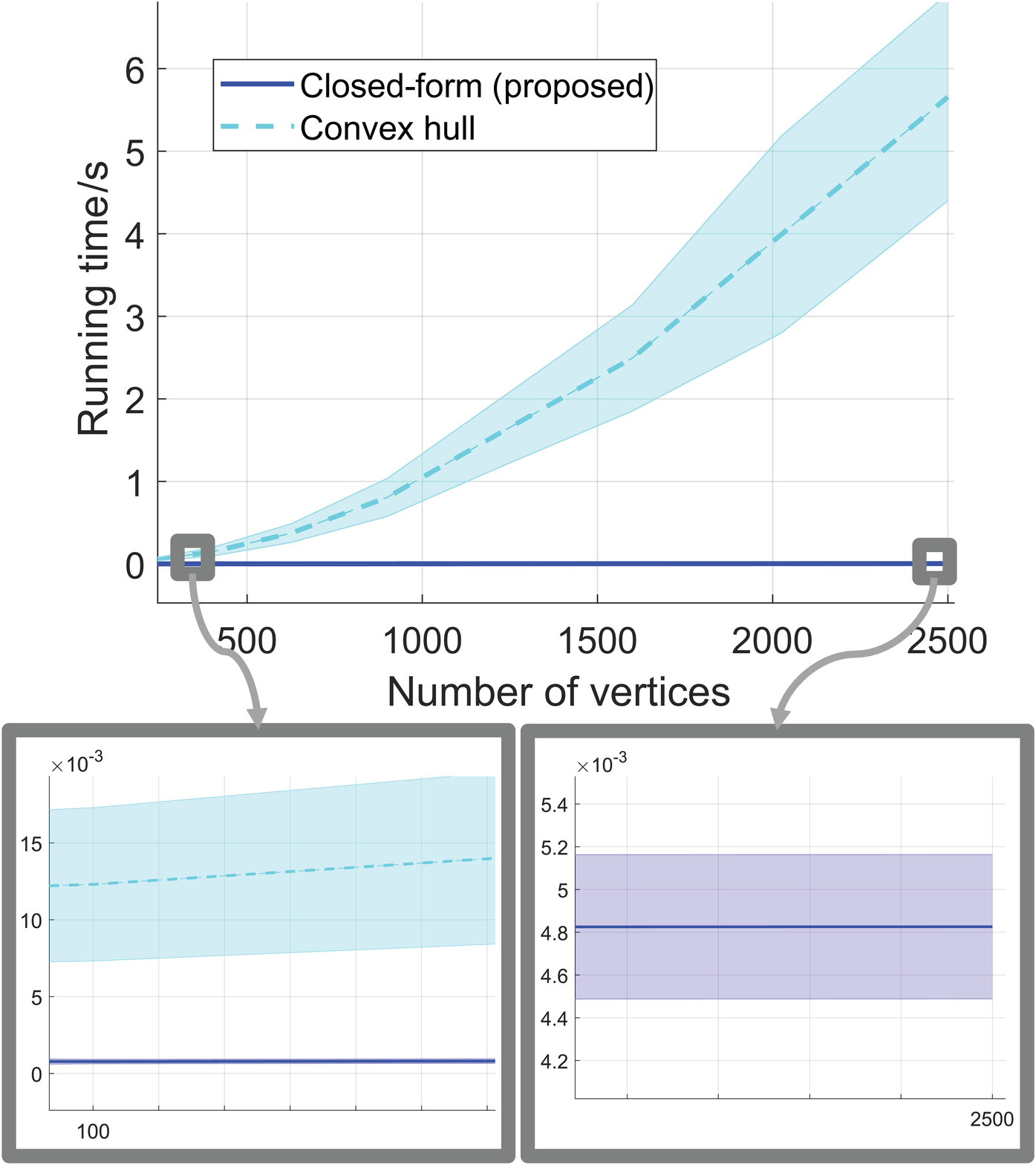}}
\caption{Running time comparisons for different algorithms in computing Minkowski sums between two superquadrics (SQ case).}
\label{fig:verification:run_time:sq}
\end{figure*}

Table \ref{tab:verification:algorithm} summarizes the applicable algorithms in different benchmark cases. And Figs. \ref{fig:verification:run_time:ee}, \ref{fig:verification:run_time:es} and \ref{fig:verification:run_time:sq} compare the running time using different algorithms in computing Minkowski sums for EE, ES and SQ cases respectively. At each discretization level, the simulations run 100 trials. 

\begin{table}[!t]
\centering
\caption{Summary of algorithms in each benchmark case}
\begin{tabular}{cccccc}
\hline
\multirow{2}{*}{Algorithm} & \multirow{2}{*}{Dim.} & Exact for & \multicolumn{3}{c}{Case} \\
\cline{4-6}
 &  & smooth surfaces? & EE & ES & SQ \\
\hline
Convex hull \citep{lozano1990spatial} & Any & & \checkmark & \checkmark & \checkmark \\
Edge sorting \citep{mark2008computational} & 2D & & \checkmark & \checkmark & \checkmark \\
Closed-form geometric \citep{yan2015closed,ruan2018path} & Any & \checkmark & \checkmark & \checkmark & \\
\textbf{Closed-form (proposed)} & Any & \checkmark & \checkmark & \checkmark & \checkmark \\
\hline
\end{tabular}
\label{tab:verification:algorithm}
\end{table}

From the comparison results, our proposed closed-form expressions for Minkowski sums achieve a linear time complexity with respect to the number of vertices on one body surface. This outperforms the method that computes the discrete Minkowski sums using convex hull, which increases much faster when more vertices are sampled. Note that an optimized convex hull algorithm for point sets is expected to help improving the performance of \emph{convex hull} method. But since this benchmark simulations only intend to provide a general impression of the performance for our proposed closed-form solution, a more detailed comparison for different convex hull algorithms is beyond the scope of this article. On the other hand, \emph{edge sorting} algorithm in 2D cases achieves linear complexity with respect to the number of vertices. Ours is faster for the simpler elliptical shapes, but slower in the case of two superellipses. This is mainly due to the complexity of computing $\tilde{\tilde{\ff}}_2(\mm_1)$, which involves calculating the mappings between gradient $\mm$ and hypersphere $\uu$. But again, both \emph{convex hull} and \emph{edge sorting} algorithms do not provide exact results for the objective studied in this work. This is because both smooth surfaces are approximated as inscribed polygons or polyhedra defined by sampled points on the surfaces. In the case that involving ellipsoids, this work achieves very similar performance with \emph{closed-form geometric}, especially in the 3D cases. However, this work is more general and not limited to ellipsoidal bodies. Furthermore, the number of the points on Minkowski sums boundary only depends on one body where the normal/gradient is computed. This explains the linear complexity of the proposed method with respect to the number of vertices on one body surface. Through these comparisons, the expression proposed in this work can compute the Minkowski sums between two convex bodies with smooth surfaces in a very efficient way, which provides quite an advantage especially in the 3D case. 

\section{Applications} \label{sec:application}
This section studies two potential applications that our proposed closed-form Minkowski sums can be used. Firstly, the generations of configuration space obstacles (C-obstacles) in robot motion planning problems are demonstrated in both planar and spatial cases. Then, the potential improvements for the queries of contact status and distance between two superquadrics are shown.

\subsection{Point-based configuration space obstacles generations}
One of the important applications of Minkowski sums is to generate configuration space obstacles \citep{lozano1990spatial} in robot motion planning algorithms. The idea is to shrink the robot into a single point and inflate all obstacles by the Minkowski sums boundaries. Then a motion planner can be developed for the point robot that avoids the inflated obstacles, which are called \emph{configuration-space obstacles (C-obstacles)}. If the robot is bounded by a circle or sphere, then, the current space is the configuration space of the robot. But if the robot is enclosed by other kinds of geometric shapes, i.e. superquadrics, the orientation matters. In this case, multiple \emph{slices} \citep{lien2008hybrid} of the configuration space are required to be generated, each of which requires the computations of the C-obstacles using Minkowski sums. Such a class of planners are effective when the robot is rigid, i.e. mobile robot, drone, underwater vehicle, etc, where the Minkowski sums computation take a very essential role.

\begin{figure*}
\centering
\subfloat[Obstacles in workspace.]{\includegraphics[scale=0.4]{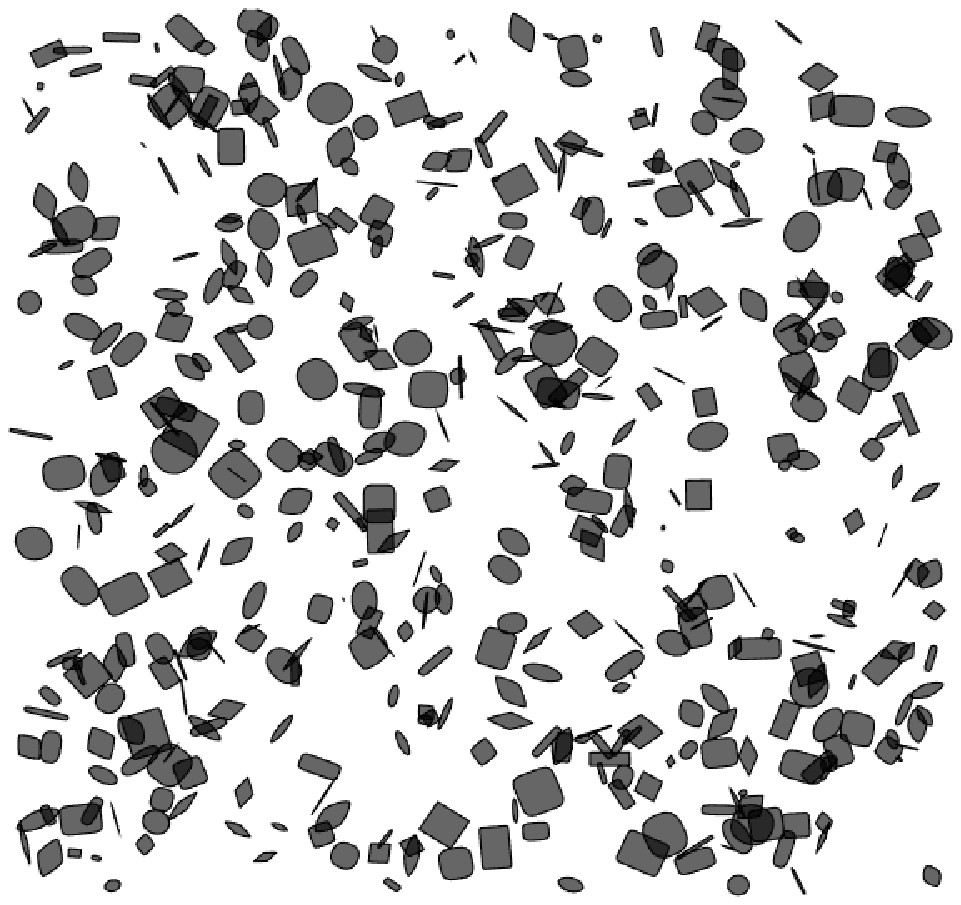}}
\subfloat[C-obstacles in different slices.]{\includegraphics[scale=0.4]{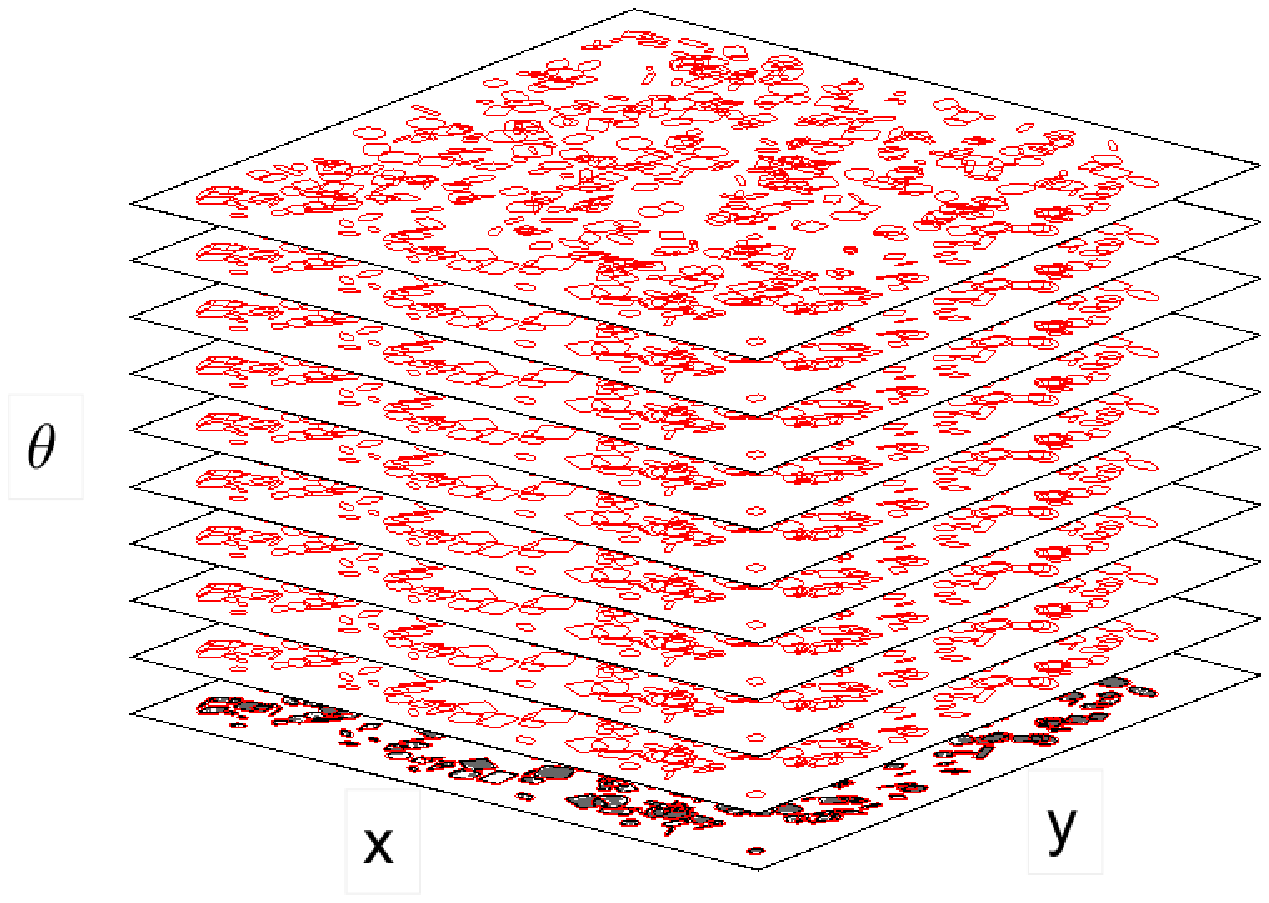}}
\caption{Demonstration of the C-obstacles generation in SE(2). The figures show the trial when 500 obstacles are located arbitrarily in the environment and the corresponding C-obstacles for 10 sampled orientations of the robot. Obstacles in the workspace are represented by black superellipses, and the calculated Minkowski sums in each C-slice are represented in red.}
\label{fig:app:c_obstacle_generation}
\end{figure*}

Here we show that using our proposed closed-form Minkowski sums, C-obstacles for rigid-body robots in multiple orientations can be calculated efficiently in both 2D and 3D cases. Both the robot and obstacles are generated as superquadric bodies with arbitrary shape parameters. We measure the computational time for generating the C-obstacle boundary points for a set of orientations of the robot and multiple numbers of obstacles. The Minkowski sums for all obstacles are computed in pairwise with the robot in a loop. Figure \ref{fig:app:c_obstacle_generation} demonstrates the simulation for C-obstacles generation in SE(2) configuration space. This demonstration includes 500 obstacles in the environment with arbitrary shapes and poses. And 10 slices of the configuration space with the computed C-obstacles are shown.

Table \ref{tab:app:c_obstacle_generation} summarizes the parameters and the running time results of generating C-obstacles in both 2D and 3D cases. The parameters include the numbers of obstacles, computed Minkowski sum boundary points and the orientations of the robot. And for the running time, the total computational time for all the points generated on all obstacles and the averaged time for each point on the Minkowski sum boundary are measured. For the 2D case, orientations are discretized as 50 heading angles of the robot within the range $[-\pi, \pi]$. And for the 3D case, orientations are sampled from the spatial rotation group, i.e. $\SO(3)$. We choose an almost uniform deterministic sampling strategy using double-coset decomposition for this discretization, which is recently proposed in \citep{wulker2019quantizing}. There are totally 3600 elements to discretize the orientation of the robot \footnote{The code to generate these 3600 samples can be referred to \url{https://github.com/ruansp/quantization-double-coset}.}.

\begin{table}[!t]
\centering
\caption{Parameters and running time for C-obstacle generations in both 2D and 3D workspaces.}
\begin{tabular}{cccccc}
\hline
Dim. & Num. & Num. & Num. & Total & Avg. time \\
& obstacles & points & orientations & time (s) & per point ($\mu s$) \\
\hline
2D & 50 & 50 & 50 & 0.36 & 2.860 \\
2D & 1000 & 50 & 50 & 6.41 & 2.563 \\
2D & 10000 & 50 & 50 & 64.26 & 2.570 \\
3D & 50 & 100 & 3600 & 95.54 & 5.308 \\
3D & 200 & 100 & 3600 & 375.02 & 5.209 \\
3D & 1000 & 100 & 3600 & 1907.85 & 5.299 \\
\hline
\end{tabular}
\label{tab:app:c_obstacle_generation}
\end{table}

The simulation shows a potential application of our closed-form Minkowski sums computations in generating C-obstacles, which is essential in developing efficient motion planning algorithms for rigid-body robots. The performance are based on a Matlab implementation without GPU or parallel computing accelerations, which achieves a level of seconds for hundreds of obstacles and minutes for thousands of obstacles. On average, the closed-form computation for each boundary point is in the level of microseconds. With these points on the boundary of C-obstacles, the collision-free configuration space can then be characterized, which leads to various efficient motion planners \citep{lien2008hybrid,ruan2018path}.

\subsection{Collision detection between two superquadrics}
Another potential application of using the proposed closed-form Minkowski sums is collision detection. In this subsection, we demonstrate possible ways of querying the contact status and proximity distance between two convex smooth bodies. Superquadrics are used as examples.

\begin{figure*}[!t]
\centering
\subfloat[Common Normal]{\includegraphics[scale=0.5, trim=70 70 60 60, clip]{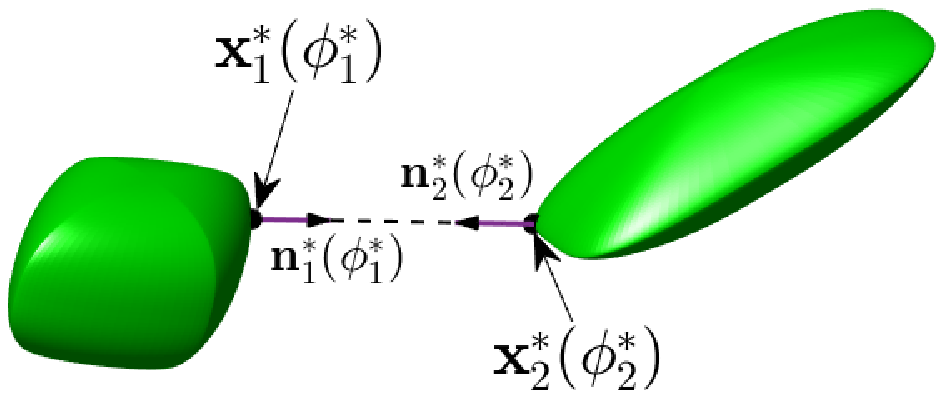}
\label{fig:app:collision:demo:cn}}
\subfloat[Minkowski (Ray)]{\includegraphics[scale=0.5, trim=50 70 60 60, clip]{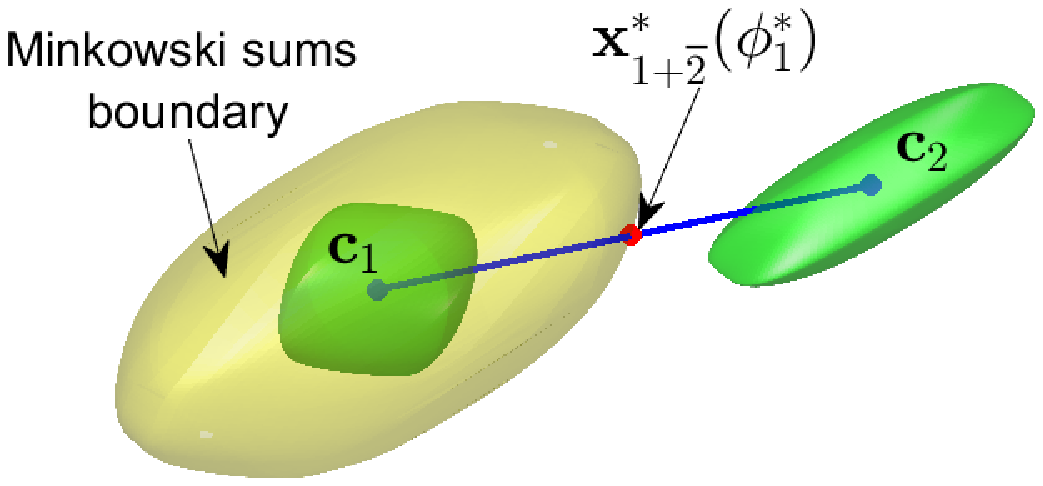}
\label{fig:app:collision:demo:mink_ray}}
\subfloat[Minkowski (Normal)]{\includegraphics[scale=0.5, trim=45 70 60 60, clip]{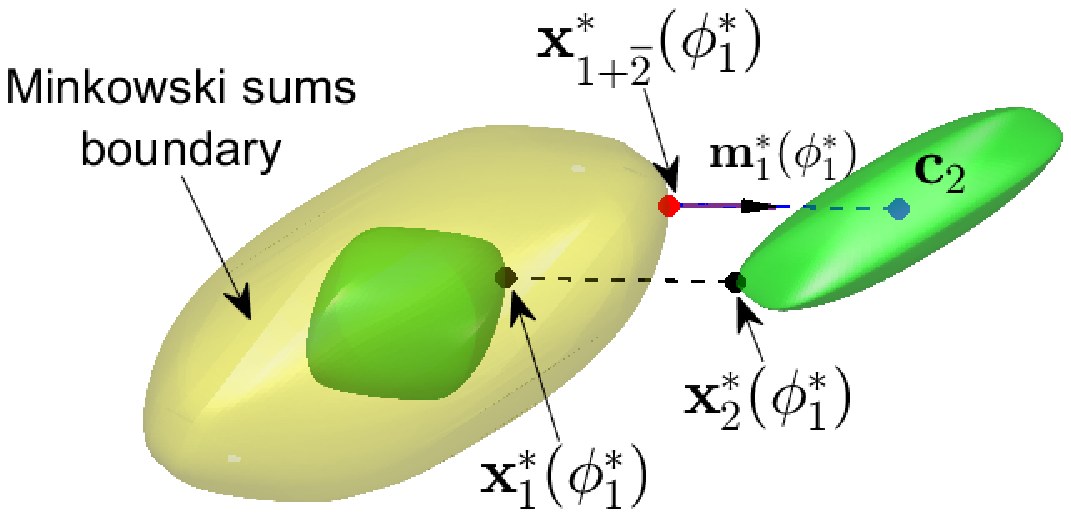}
\label{fig:app:collision:demo:mink_normal}}
\caption{Geometric interpretations of different optimization-based collision detection algorithms. The geometric bodies are in green color and Minkowski sums boundary is in transparant yellow color. Superscript ``$*$'' is labelled to the optimal solution for variables, surface points normal/gradient vectors.}
\label{fig:app:collision:demo}
\end{figure*}

\subsubsection{Review of some existing algorithms}
There are many elegant collision detection algorithms designed for superquadrics, such as using implicit surfaces with convex optimization \citep{bloomenthal1997introduction,chakraborty2008proximity} and common normal concept \citep{gonccalves2017benchmark}. The variables to be optimized are the parameterized points on both bodies. For example, given the implicit functions of $B_1$ and $B_2$, i.e., $\Psi_1(\cdot)$ and $\Psi_2(\cdot)$ respectively, the witness points that have minimum distance can be obtained by solving
\begin{equation}
\begin{aligned}
\min_{\xx_1, \xx_2} \,\, & \| \xx_1 - \xx_2 \|_2^2 \\
\text{subject to} \,\, & \Psi_1(\xx'_1) \leq 0 \,, \\
& \Psi_2(\xx'_2) \leq 0 \,,
\end{aligned}
\label{eq:app:collision:implicit}
\end{equation}
where $\xx_1, \xx_2$ are the points expressed in the global frame and $\xx'_1, \xx'_2$ are the corresponding points in the local frame of body $B_1$ and $B_2$ respectively. This is a convex optimization problem and can be solved by the interior-point algorithm \citep{chakraborty2008proximity}. The two bodies overlap when the optimal objective function equals zero and separated otherwise. Furthermore, a necessary condition of obtaining the minimum distance between two bodies is based on the common normal concept \citep{lopes2010mathematical, gonccalves2017benchmark}. Efficient algorithms have been proposed by solving
\begin{equation}
\begin{pmatrix}
{\bf n}_1(\pphi_1) \times {\bf n}_2(\pphi_2) \\
{\bf n}_1(\pphi_1) \times [\xx_2(\pphi_2) - \xx_1(\pphi_1)]
\end{pmatrix}
= {\bf 0} \,,
\label{eq:app:collision:common_normal}
\end{equation}
where ${\bf n}_1, {\bf n}_2$ are outward-pointing normal vector at surface points $\xx_1, \xx_2$ respectively. The surface points and normal vectors are viewed in the global frame and parameterized by their corresponding spherical coordinates $\pphi_1, \pphi_2$. This nonlinear equation can also be solved by numerical methods, e.g. Levenberg-Marquardt (L-M) algorithm. The optimal solution of an example trial is demonstrated in Fig. \ref{fig:app:collision:demo:cn}. The length of the line segment connecting two witness points (dashed line) is the solved minimum Euclidean distance between the two bodies.

\subsubsection{Applications of our proposed closed-form Minkowski sums on collision detection problems}
The collision detection can be equivalently and efficiently solved by applying our proposed closed-form Minkowski sums expressions. Previous work has shown the superiority for querying contact status between an ellipsoid and a superquadric body using closed-form Minkowski sums \citep{ruan2019efficient}. This method uses the spirit of ray casting. The same algorithmic pipeline can be used when combining with the expression in this article. The idea is to solve for $\pphi_1$ in the nonlinear equation
\begin{equation}
\xx_{1+\overline{2}}(\pphi_1) \times ( {\bf c}_2 - {\bf c}_1 ) = {\bf 0},
\label{eq:app:collision:mink_ray}
\end{equation}
where ${\bf c}_1, {\bf c}_2$ are the center points of $B_1, B_2$ respectively. This method can return the contact status between the two bodies. Figure \ref{fig:app:collision:demo:mink_ray} geometrically demonstrates the idea. The solved point $\xx^*_{1+\overline{2}}$ is the intersection between the Minkowski sums boundary and the line segment connecting the centers of two bodies. 

To solve for the Euclidean distance between two bodies, we propose another objective function by using the ``common normal'' concept. Since the Minkowski sums boundary guarantees anti-parallelism of normal/gradient vectors, Eq. \eqref{eq:app:collision:common_normal} can be simplified into solving for $\pphi_1$ from
\begin{equation}
{\bf m}_1(\pphi_1) \times [{\bf c}_2 - \xx_{1+\overline{2}}(\pphi_1)] = {\bf 0} \,.
\label{eq:app:collision:mink_normal}
\end{equation}
Note that one can also use the normalized ${\bf n}_1$ and Eq. \eqref{eq:mink_sum_norm} in the objective function. But here, we demonstrate the idea using the more general un-normalized gradient $\mm_1$ and Eq. \eqref{eq:mink_sum_non_norm} for $\xx_{1+\overline{2}}$. The geometric interpretation of this algorithm is demonstrated in Fig. \ref{fig:app:collision:demo:mink_normal}. The solved gradient $\mm^*_1$ is co-linear with the line connecting the center of $B_2$ and the optimal point $\xx^*_{1+\overline{2}}$ on Minkowski sums boundary. The resulting withness points $\xx^*_1$ and $\xx^*_2$ are identical with those solved by the common normal method, i.e. Eq. \eqref{eq:app:collision:common_normal}.

\subsubsection{Comparisons on solving collision detection queries}
The performance of solving the aforementioned nonlinear equations are compared quantitatively. The algorithms to solve Eqs. \eqref{eq:app:collision:implicit}, \eqref{eq:app:collision:common_normal}, \eqref{eq:app:collision:mink_ray} and \eqref{eq:app:collision:mink_normal} are denoted as \emph{Implicit}, \emph{Common Normal}, \emph{MinkSum (Ray)} and \emph{MinkSum (Normal)}, respectively. To make the comparison as fair as possible, apart from Eq. \eqref{eq:app:collision:implicit}, all the other equations are solved by L-M algorithm. These optimization-based methods use exact expressions of a superquadric model. For the special case of ellipsoids, an algebraic condition of separation (denoted as \emph{ASC}) has been proposed to solve for the collision status only \citep{wang2001algebraic}. And for convex polyhedra, \emph{GJK} is a very popular and effective algorithm to solve for both collision status and distance \citep{gilbert1988fast}. Although not using the same type of computational framework, to make the comparisons more comprehensive, these two methods are also included in the benchmark studies\footnote{GJK implementation is retrieved from \url{https://github.com/mws262/MATLAB-GJK-Collision-Detection}, and ASC implementation is based on \url{https://www.mathworks.com/matlabcentral/fileexchange/32172-are-two-ellipsoids-in-contact-algebraic-separation-condition-for-ellipsoids}.}. The convex polyhedra for GJK method are generated using discrete surface points of the superquadrics, the same as in Sec. \ref{sec:numerical_verifications:poly_sq}. Based on the previous simulations, to fit the superquadrics well, the number of vertices for each convex polyhedron is chosen to be 100. In the benchmark simulations, collision detection between two superquadrics at random orientations and positions are solved. For each method, a total of $10^4$ random poses for the two bodies are conducted. The averaged solving time among all trials for different algorithms are summarized in Tab. \ref{tab:app:collision:time}.

\begin{table}[!t]
\centering
\caption{Running time comparisons among collision detection algorithms for different pairs of objects}
\begin{tabular}{ccc}
\hline
Algorithm & Ellipsoids (ms) & Superquadrics (ms) \\
\hline
GJK & 18.4 & 18.3 \\
ASC & {\bf 0.133} & -- \\
Implicit & 27.0 & 42.3 \\
Common Normal & 4.1 & 7.1 \\
\textbf{MinkSum (Ray)} & 3.5 & 5.2 \\
\textbf{MinkSum (Normal)} & 3.5 & {\bf 5.1} \\
\hline
\end{tabular}
\label{tab:app:collision:time}
\end{table}

All the collision detection queries for all algorithms can be solved in the level of milliseconds. The \emph{Implicit} method runs about 5-10 times slower than the others, which can be reasoned as the different performances of optimization solvers. But the better performance of both \emph{MinkSum} methods compared to the \emph{Common Normal} method can be mainly due to the usage of the closed-form expression of Minkowski sums. The variables as well as the dimension of the objective function to be optimized are reduced into the half of that in the \emph{Common Normal} method. In fact, the common normal concept is implicitly encoded in the closed-form Minkowski sums expression. The \emph{GJK} algorithm uses the discrete convex polyhedron approximation of the superquadric model. Its performance highly depends on the number of vertices of the polyhedron. To balance the efficiency and approximation errors, the number is chosen as small as possible such that the relative volume error is less that $10\%$, as Sec. \ref{sec:numerical_verifications:poly_sq} studied. It is shown that with this discretization, GJK runs faster than the \emph{Implicit} method but slower than \emph{Common Normal} method and methods using our closed-form Minkowski sums expressions. \emph{ASC} method provides a very efficient way to determine the separation condition of two ellipsoids, but it cannot be generalized to other geometric models like other algorithms do. For the case of two ellipsoids, the \emph{MinkSum} methods still outperform the others except for \emph{ASC}. Thus, the efficiency of computing the closed-form Minkowski sums can potentially help speeding up the optimization-based collision detection algorithms.

\section{Conclusion and Future Work} \label{sec:conclusion}

This article introduces a novel exact and closed-form parameterization of Minkowski sums between two positively curved bodies in $d$-dimensional Euclidean space. The boundary surface of each body is parameterized as a function of the surface gradient. It is shown that the Minkowski sums can be derived in closed-form based on the surface gradient that can be either normalized (i.e. becomes the outward normal vector) or un-normalized. The results proposed in this article, both in canonical form and with linear transformations, are identical with those using geometrical interpretation in the case of two ellipsoids. And in general, these expressions can be applied to pairs of bodies that are enclosed by general smooth surfaces with all sectional curvatures being positive at every point. Numerical simulations are conducted in the case of two superquadrics, showing the correctness and efficiency of the proposed closed-form expression. Applications in generating configuration-space obstacles for motion planning algorithms and improving optimization-based collision detection algorithms are introduced and demonstrated numerically.

The current work still has some limitations. For example, we only derived the closed-form expression for bodies under linear transformations, but the possibility of applying nonlinear deformations (e.g. tapering) is unknown. If no exact expression can be obtained, the upper or lower bounds might be derived. Also, the extension to non-convex bodies enclosed by smooth surfaces still remains for further investigation. In this case, the inverse Gauss map from surface normal to surface point is not unique, so a direct extension of the current framework might not be possible. But it might be able to use the idea from convex decomposition-based methods to address the non-convexity of the bodies.

In the future, the authors seek to solve the aforementioned limitations, i.e. exploring possibilities of applying general deformations to both bodies and extending to non-convex bodies. We would also like to dig more into the possible improvements of applying the closed-form Minkowski sums in various collision detection algorithms.

\section*{Acknowledgement}
This work was performed under National University of Singapore Startup grants R-265-000-665-133 and R-265-000-665-731 and National University of Singapore Faculty Board funds C-265-000-071-001. When the authors were at Johns Hopkins University, they were supported by U.S National Science Foundation grant IIS-1619050. The ideas expressed in this article are solely those of the authors.

\bibliographystyle{elsarticle-num}
\bibliography{reference}

\end{document}